\documentclass[11pt]{article}
\usepackage{dsfont}
\usepackage{pxfonts}
\usepackage{yfonts}
\usepackage{graphicx}
\usepackage{relsize}

\def\bel{\begin{equation}\label}
\def\eeq{\end{equation}}
\def\ds{\displaystyle}
\def\endproof{\hphantom{MM}
\hfill\llap{$\square$}\goodbreak}

\def\mt{\longrightarrow}
\def\v{\vskip 1em}

\def\R{\mathbb{R}}

\def\Cx{\mathbb{C}}

\def\exp{{\bf exp}}
\def\Re{{\bf Re}}
\def\Im {{\bf Im}}

\def\D{{\bf D}}

\def\A{{\bf A}}

\def\p{{\partial}}
\def\a{{\bf a}}
\def\b{{\bf b}}

\def\i{{\bf i}}
\def\Tilde{\widetilde}
\def\Hat{\widehat}

\def\I{{\bf I}}
\def\II{{\bf II}}
\def\III{{\bf III}}
\def\IV{{\bf IV}}
\def\M{{\bf M}}
\def\V{{\bf V}}
\def\U{{\bf U}}

\def\Cup{{\bigcup}}

\def\xi{{\xiup}}
\def\eta{{\etaup}}
\def\tau{{\tauup}}
\def\alpha{\alphaup}
\def\lambda{{\lambdaup}}
\def\omega{\omegaup}
\def\varphi{{\varphiup}}

\def\adj{\hbox{adj}}
\def\exp{\hbox{\bf exp}}

\newtheorem{remark}{Remark}[section]

\begin{document}
\title{\bf Generalized Dirichlet to Neumann Maps for\\ 
Linear Dispersive Equations on Half-Line}
\v
\author{Athanassios S. Fokas$^{*}$ and Zipeng Wang$^\dagger$
 \\ \\ 
($\ast$)~Department of Applied Mathematics and Theoretical Physics, \\
University of Cambridge, Cambridge, CB3 0WA, United Kingdom.
\\\\
e-mail: T.Fokas@damtp.cam.ac.uk
\\\\
($\dagger$)~Cambridge Centre for Analysis, University of Cambridge,\\ 
Cambridge, CB3 0WA, United Kingdom.
\\\\
e-mail:  z.wang@maths.cam.ac.uk}
\maketitle
\begin{abstract}
A large class of initial-boundary problems of linear evolution partial differential equations 
formulated on the half-line is analyzed via the unified transform method. In particular, explicit representations are presented for the generalized Dirichlet to Neumann maps. Namely, the determination of the unknown boundary values when an essential set of initial and boundary data is given.
\end{abstract}

\section{Introduction}
\setcounter{equation}{0}
In this paper, we study a class of initial-boundary value problems of linear evolution equations formulated on the half line, by using so-called the {\it unified transform} method introduced in \cite{ASF1}. A major difficulty of solving  certain initial-boundary value problems  stems from the fact that the solution representation requires all boundary values, whereas only a subset of them is prescribed as boundary conditions. The determination of the unknown boundary values in terms of given data is often called the {\it Generalized Dirichlet to Neumann maps}.

In particular, we consider the P.D.E
\bel{PDE}
\Big(\p_t+\omega(-\i\p_x)\Big)q(t,x)~=~0,\qquad 0<t<T,\qquad 0<x<\infty
\eeq
where  $\omega$ is defined by
\bel{omega xi}
\omega(\xi)~=~a_n \xi^n+a_{n-1}\xi^{n-1}+\cdots+a_1\xi+a_o
\eeq
for which $\{a_j\}_{j=0}^n$ are complex constant coefficients.  

We assume that the initial condition
\bel{initial}
q_o(x)~\doteq~q(0,x)
\eeq
has sufficient decaying as $x\mt\infty$, as discussied in \cite{ASF1}-\cite{ASF}.
Here, we analyze the generalized Dirichlet to Neumann maps for the following initial-boundary value problems: 

({\bf a}) Let 
\bel{omega*}
\omega(\xi)~=~a_n\xi^n
\eeq
where a subset of $n-N$ boundary values is prescribed as boundary conditions. The exact value of $N$ will be specified later. 
We set 
\bel{U-V}
U\cup V~=~\Big\{0,1,\ldots,n-1\Big\}
\eeq
where
\bel{UV}
U=\Big\{ u_1,u_2,\ldots,u_{n-N}\Big\},\qquad V=\Big\{ v_1,v_2,\ldots,v_N\Big\}.
\eeq
 Let
\bel{given conditions}
g_l(t)~\doteq~\p_x^{u_l}q(t,0),\qquad u_l\in U 
\eeq
which are given smooth functions, and are compatible with $q_o(x)$ at $x=0$. Namely, 
\bel{g_l(0)}
g_l(0)~=~q_o(0),\qquad u_l\in U.
\eeq
 We will present explicit formulae for the remaining unknown boundary values 
\bel{boundary values}
\p_x^{v_j}q(0,t),\qquad v_j\in V
\eeq
in terms of $\Big\{g_l(t)\Big\}_{l=1}^{n-N}$. Some examples were computed in \cite{ASF} for $n=2,3$.

({\bf b}) Let $\omega$ in the general form (\ref{omega xi}) for $n\leq5$ but assume the {\it canonical} boundary conditions. 
Namely, $q(t,0)$ and its first $n-N-1$ derivatives are prescribed as boundary conditions:   
 \bel{given conditions canonical}
 \p_x^{u_l} q(t,0)~=~g_l(t),\qquad u_l\in \Big\{0,1,2,\ldots,n-N-1\Big\}.
 \eeq
 We will present explicit formulae for the unknown boundary values in terms of $\Big\{g_l(t)\Big\}_{l=1}^{n-N}$.
 These formulae involve the solution of an algebraic equation of order $n-1$, which in general can be solved in terms of radicals only when $n\leq5$.    
 
The unified transform ( or Fokas transform) was introduced in  \cite{ASF1}, \cite{ASF2} ( see also the book \cite{ASF} and the reviews \cite{F.S}, \cite{D.T.V} ). 
The implementation of the unified transform to evolution equations on the half-line and the finite interval is discussed in \cite{ASF3}, \cite{M.F} and \cite{F.P}-\cite{P1} respectively. 
The case of periodic initial condition is discussed in \cite{T.D}. The large $t$ asymptotics of evolution PDEs on the half-line is analyzed in \cite{F.S1}-\cite{D}. 
The numerical implementation of the unified transform to evolution PDEs is discussed in \cite{F.F} and \cite{V-C}. 
Rigorous results are presented in \cite{FS}. The implementation of the unified transform to evolution PDEs in the two dimensional space is discussed in \cite{ASF3}, \cite{F.K} and \cite{M.F1}.  
Implications of the unified transform in the area of spectral theory are discussed in \cite{ASF4}-\cite{S}. The implementation of the unified transform to several problems of physical significance is 
discussed in \cite{F.S2}-\cite{A.S.P.S}. Systems of evolution PDEs are considered in \cite{Spen} and \cite{Dima}.

\section{Formulation on the Main Results}
\setcounter{equation}{0}
 In this section, we introduce our main results after certain preliminary settings, for which more background discussions can be found in the book \cite{ASF}. 

Observe that equation (\ref{PDE}) admits the family of explicit solutions: 
\bel{solution}
\exp\Big\{\i x\xi-\omega(\xi)t\Big\},\qquad0<x<\infty,\qquad0<t<T
\eeq
where $\xi\in\Cx$ . The convergence of the solution is provided by $\xi\in\Cx^+\cup\R$ and $\Re{\omega(\xi)}\ge0$. Taking into consideration that $\omega(\xi)$ asymptotes to $a_n\xi^n$ as $|\xi|\mt\infty$, we thus assume that $\Re {a_n}\ge0$ if $n$ is even, and $\Re{a_n}=0$ if $n$ is odd. Without losing of generality, we let $|a_n|=1$. 

Let $\D$ to be the {\it principal domain}:
\bel{principal domain}
\D~=~\Big\{\xi\in\Cx~\colon~\Re\Big(a_n\xi^n\Big)<0\Big\}
\eeq
which is an union of $n$ sectors in $\Cx$. 

Throughout the rest of the paper, we fix the notations 
\bel{a_n,xi}
 a_n=\cos\varphi+\i\sin\varphi\qquad\hbox{and}\qquad \xi=R(\cos\vartheta+\i\sin\vartheta).
 \eeq
Notice that $\Re{a_n\xi^n}<0$ implies $\cos(\varphi+n\vartheta)<0$. 
The $n$ sectors of $\D$ are characterized by  
\bel{Degrees}
\begin{array}{lr}\ds
\Bigg( \vartheta+{\varphi-2m\pi\over n}\Bigg)~\in~\Bigg({\pi\over 2n}~,~{3\pi\over 2n}\Bigg)\qquad m=0,1,\ldots,n-1.
\end{array}
 \eeq
 Moreover, $\Re{a_n}\ge0$ if $n$ is even and $\Re{a_n}=0$ if $n$ is odd together imply that
\bel{varphi}
\varphi~\in~\left[-{\pi\over 2}~,~{\pi\over 2}\right].
\eeq
A direct computation shows that there are exactly $N$ sectors of $\D$ lying in $\Cx^-$, whereas
\bel{N}
N~=~\left\{\begin{array}{lr}
\ds {n\over 2}\qquad\qquad \hbox{if}~n~\hbox{is even},
\\\\
\ds {n-1\over 2}\qquad~ \hbox{if}~n~\hbox{is odd}~~a_n=\i,
\\\\
\ds {n+1\over 2}\qquad~ \hbox{if}~n~\hbox{is odd}~~a_n=-\i.
\end{array}\right.
\eeq
The solution $q(t,x)$ of the initial boundary value problems ($\a$) and ($\b$) is well-poseded 
provided that there are essentially $N$ boundary values given. See \cite{ASF} and \cite{FS} for references. 

We write
\bel{prin domain +/-}
 \D^-~=~\D\cap\Cx^-~=~\Cup_{j=1}^N ~\D^-_j\qquad\hbox{and}\qquad
\D^+~=~\D\cap\Cx^+~=~\Cup_{k=1}^{n-N} ~\D^+_k
\eeq
where $\D^+_k$ and $\D^-_j$ represent the sectors as subdomains of $\D\cap\Cx^+$ and $\D\cap\Cx^-$ respectively. These sectors are pairwisely disjoint. We number them consecutively by starting with the one with the smallest degree. 

We write the boundary of $\D^+$ by $\p \D^+$, whereas
\bel{prin domain boundary}
\p\D^+~=~\Cup_{k=1}^{n-N} \p\D^+_k
\eeq
with a counterclockwise orientation. 
Each $\p\D^+_k$ consists of the rays $\mathcal{R}^1_k$ and $\mathcal{R}^2_k$ with $\mathcal{R}^1_k$ emanating from the origin and $\mathcal{R}_k^2$ directed towards the origin. 
Geometrically, we always have $\mathcal{R}_k^1$ on the right of $\mathcal{R}_k^2$ for every $k=1,2,\ldots,n-N$. 

By using the notations
\bel{vartheta}
\vartheta^1_k~=~\arg\mathcal{R}_k^1\qquad\hbox{and}\qquad \vartheta^2_k~=~\arg\mathcal{R}^2_k,
\eeq
Equations (\ref{Degrees})-(\ref{varphi}) imply respectively that
\bel{Ray degree}
\vartheta_k^1~=~{\pi\over 2n}-{\varphi-2(k-1)\pi\over n},\qquad \vartheta_k^2~=~{3\pi\over 2n}-{\varphi-2(k-1)\pi\over n}
\eeq
for $k=1,2,\ldots,n-N$. In particular, we have
\bel{vartheta Est}
\vartheta_1^1~=~\left\{\begin{array}{lr}\ds
{\pi-2\varphi\over 2n}\qquad n~\hbox{is even},
\\\\ \ds
{\pi\over n}\qquad n~\hbox{is odd},~~a_n=-\i
\\\\ \ds
0\qquad n~\hbox{is odd},~~a_n=\i
\end{array}\right.
\qquad\qquad
\vartheta_{n-N}^2~=~\left\{\begin{array}{lr}\ds
\pi-{\pi-2\varphi\over 2n}\qquad n~\hbox{is even},
\\\\ \ds
\pi-{\pi\over n}\qquad n~\hbox{is odd},~~a_n=-\i
\\\\ \ds
\pi\qquad n~\hbox{is odd},~~a_n=\i.
\end{array}\right.
\eeq
Let 
\bel{rho}
\rho~=~\exp\left({2\pi\i\over n}\right).
\eeq
Define the $N\times N$ Vandermonde matrix
\bel{V}
\V(\rho)=\left[\begin{array}{lr}
\rho ^{n-v_1-1}~~~~~~~\rho ^{n-v_2-1}~~~~\cdots~~~\rho ^{n-v_N-1}
\\\\
\rho^{2(n-v_1-1)}~~~~\rho^{2(n-v_2-1)}~\cdots~~\rho^{2(n-v_N-1)}
\\\\
~~~~~~~~~~~~~~~~~\vdots ~\qquad\qquad~\ddots
\\\\
\rho^{N(n-v_1-1)}~~~\rho^{N(n-v_2-1)}~\cdots~ \rho^{N(n-v_N-1)}
\end{array}\right]
\eeq
whose determinant is nonzero. 

Let $\V_{ji}$ to be the $(i,j)$-th principal minor of $\V$. 
When the matrix is a scalar, we take its $(1,1)$-th principal minor to be $1$.
From (\ref{V}) we have
\bel{det V Cremer rule}
\det\V(\rho)~=~\sum_{i=1}^N (-1)^{j+i} \det\V_{ji}(\rho)~ \rho^{i(n-v_j-1)},\qquad j=1,2,\ldots,N.
\eeq
Let $\V^{jl}$  to be the matrix obtained from the matrix $\V$, after replacing  its $j$-th column by the column vector 
$\left(\rho^{i(n-u_l-1)}\right)_{1\leq i\leq N}$. 

Our first main result regarding problem ($\a$) is given below.

{\bf Theorem One:} {\it Let $q(t,x)$ satisfy the PDE in (\ref{PDE}) with $\omega(\xi)$ defined in (\ref{omega*}). Given
\bel{given condition 1}
q(0,x)~=~q_o(x),~~~0<x<\infty; \qquad \p_x^{u_l}q(t,0)~=~g_l(t), ~~~0<t<T,~~u_l\in U
\eeq
where $q_o(x)$ and $g_l(t)$ are defined in (\ref{initial}) and (\ref{given conditions}), 
the unknown  boundary values can be determined by the following formulae:
\bel{Result1} 
\begin{array}{lr}\ds
 2\pi\p^{v_j}_x q(t,0)~=~-{1\over n-N}\sum_{k=1}^{n-N}\int_{\p \D^+_k}\sum_{i=1}^N(-1)^{i+j}\left({\det \V_{ji}\over \det\V}\right)(\rho^{n-N+1-k})n\Hat{q}_o(\rho^{n-N+i-k}\xi)(\i\xi)^{v_j}
 \\\\ \ds
 ~+~{1\over n-N}\sum_{k=1}^{n-N} \int_{\p \D^+_k}\sum_{u_l<v_j}\left({\det\V^{jl}\over\det\V}\right)(\rho^{n-N+1-k})nq_o(0)(\i\xi)^{v_j-u_l-1}e^{-\omega(\xi)t}d\xi
 \\\\ \ds 
 ~+~ \sum_{v_j<u_l}\i^{v_j-u_l+1} \Lambda_{jl}\Gamma\left({v_j-u_l+n\over n}\right)\int_0^t {g_l(\tau)d\tau\over (t-\tau)^{v_j-u_l+n\over n}}
 ~-~\sum_{u_l<v_j}\i^{v_j-u_l+1}\Lambda_{jl}\Gamma\left({v_j-u_l\over n}\right) \int_0^t {\dot{g}_l(\tau)d\tau\over (t-\tau)^{v_j-u_l\over n}} 
 \end{array}
 \eeq
 for every $v_j\in V$, where $\Gamma$ is Gamma function and 
\bel{Result1'} \begin{array}{lr}\ds
\Lambda_{jl}~=~
\left({\det\V^{jl}\over\det\V}\right)(\rho^{n-N})\exp\left(\i(v_j-u_l)\left(\vartheta_1^1-{\pi\over 2n}\right)\right)-\left({\det\V^{jl}\over\det\V}\right)(\rho)\exp\left(\i(v_j-u_l)\left(\vartheta^2_{n-N}+{\pi\over 2n}\right)\right).
\end{array}
\eeq }
Let $\omega$ in the general form  (\ref{omega xi}) but assume $n\leq 5$. Consider the domain
\bel{Domain}
D~=~\Big\{\xi\in\Cx~\colon~\Re\Big(\omega(\xi)\Big)<0\Big\}.
\eeq
Since $\omega(\xi)\approx a_n\xi^n$ as $\xi$ is asymtotically large, the domain $D$ coincides with the principal domain $\D$ in (\ref{principal domain}) at infinity.
We define 
\bel{D_L union}
D_L~=~\Big\{ \xi\in D~\colon~|\xi|\ge L\Big\} \qquad
 D^-_L~=~D_L\cap\Cx^-~=~\Cup_{j=1}^N ~D^-_{L,j}\qquad
  D^+_L~=~D_L\cap\Cx^+~=~\Cup_{k=1}^{n-N} ~D^+_{L,k}
  \eeq
for some $L$ sufficiently large. 
Again, we number them consecutively by starting with the one with the smallest degree.
Moreover, we write the boundary of $D^+_L$ by $\p D^+_L$, whereas
\bel{domain boundary}
\p D^+_L~=~\Cup_{k=1}^{n-N} \p D^+_{L,k}
\eeq
with a counterclockwise orientation, for $L$ sufficiently large. 
By fundamental theorem of algebra, there are exactly
$n-1$ functions $z_i(\xi), i=1,2,\ldots,n-1$ other than $\xi$ itself, determined implicitly by
 \bel{Poly}
{\omega(z_i)-\omega(\xi)\over z_i-\xi}~=~0
\eeq
such that $\omega(z_i(\xi))=\omega(\xi)$ for every $i=1,2,\ldots,n-1$.
Let $\xi\in D^+_k$ for $k=1,2,\ldots,n-N$, there are exactly $N$ such functions 
\bel{z}
z_1^k\left(\xi\right),~z_2^k\left(\xi\right),~\ldots,~z_N^k\left(\xi\right)
\eeq
whose values lie inside $\Cx^-$ for $\xi$ asymptotically large. See \cite{FS} for more discussions.

Let $\omega_n(\xi)=\omega(\xi)$ and define inductively
\bel{omega_m}
\omega_m(\xi)~=~{1\over \xi}\Big(\omega_{m+1}(\xi)-a_{n-m-1}\Big),\qquad m=1,2,\ldots,n-1.
\eeq 
 The $N\times N$ alternant matrix 
\bel{A}
\A(k,\xi)~\doteq~\A(z_1^k(\xi),z_2^k(\xi),\ldots,z_N^k(\xi))~=~\left[\begin{array}{ccc}\ds
\omega_{n-v_1-1}(z_1^k(\xi))~~\omega_{n-v_2-1}(z_1^k(\xi))~\cdots~\omega_{n-v_N-1}(z_1^k(\xi))
\\\\
\omega_{n-v_1-1}(z_2^k(\xi))~~\omega_{n-v_2-1}(z_2^k(\xi))~\cdots~\omega_{n-v_N-1}(z_2^k(\xi))
\\
\vdots~\qquad\qquad~\ddots
\\
\omega_{n-v_1-1}(z_N^k(\xi))~~\omega_{n-v_2-1}(z_N^k(\xi))~\cdots~\omega_{n-v_N-1}(z_N^k(\xi))
\end{array}\right].
\eeq
Let $\A_{ji}$ to be the $(i,j)$-th principal minor of $\A$. From (\ref{A}) we have
\bel{CramerA}
\det\A(k,\xi)~=~\sum_{i=1}^N (-1)^{j+i} \det\A_{ji}(k,\xi)~ \omega_{n-v_j-1}(z_i^k(\xi)),\qquad j=1,2,\ldots,N.
\eeq
Let $\A^{jl}$ to be the matrix obtained from the matrix $\A$ after replacing its $j$-th column by the column vector 
$\left(\omega_{n-u_l-1}(z_i^k(\xi))\right)_{1\leq i\leq N}$. 
Define
\bel{lambda}
\lambda~+~\exp\left({\pi\i\over n-N}\right).
\eeq
Our second main result regarding problem ($\b$) is given below.

{\bf Theorem Two:} {\it Let $q(t,x)$ satisfy the PDE in (\ref{PDE}) with $\omega(\xi)$ defined in (\ref{omega xi}) and $n\leq 5$. Given 
\bel{given condition 2}
\begin{array}{lr}\ds
q(0,x)~=~q_o(x),~~~0<x<\infty; 
\\\\ \ds
 \p_x^{u_l}q(t,0)~=~g_l(t), ~~~0<t<T,~~u_l\in\Big\{0,1,\ldots,n-N-1\Big\}
 \end{array}
 \eeq
where $q_o(x)$ and $g_l(t)$ are defined in (\ref{initial}) and (\ref{given conditions}),  
the unknown  boundary values can be determined by the following formulae:
\bel{Result2}
\begin{array}{lr}\ds
2\pi(n-N)\p_x^{v_j}q(t,0)~=~ -\i^{v_j}\sum_{k=1}^{n-N}\int_{\p D^+_{L,k}} \sum_{i=1}^N (-1)^{i+j}\left({\det \A_{ji}\over \det\A}\right)\left(k,\xi\right)\omega'(\xi)\Hat{q}_o\left(z_i^k\left(\xi\right)\right)e^{-\omega(\xi)t}d\xi
\\\\ \ds
~-~\sum_{k=1}^{n-N}\int_{\p D^+_{L,k}}\sum_{l=1}^{n-N}\i^{v_j-u_l+1} 
 \left({\det\A^{jl}\over \det\A}\right)\left(k,\xi\right){\omega'(\xi)\over \omega(\xi)} 
q_o(0)e^{-\omega(\xi)t}d\xi
\\\\ \ds
~+~\sum_{k,l=1}^{n-N} \i^{v_j-u_l+1}~\hbox{\bf p.v}\int_{0}^{\infty}\left[\left({\det\A^{jl}\over \det\A}\right)\left(k,\lambda^k\xi\right)-\left({\det\A^{jl}\over \det\A}\right)\left(k,\lambda^{k=1}\xi\right)\right]{\omega'(\xi)\over \omega(\xi)}\left( \int_0^t e^{\omega(\xi)(\tau-t)}\dot{g}_l(\tau)d\tau\right)d\xi
\\\\ \ds 
~+~\sum_{k,l=1}^{n-N} \i^{v_j-u_l+1}\Big(g_l(t)-g_l(0)\Big)\left[ 2\pi\i\sum_{\omega(\zeta)=0,\arg\zeta\in\Theta_k} \left({\det\A^{jl}\over \det\A}\right)\left(k,\zeta\right)+\pi\i\sum_{\omega(\zeta)=0,\arg\zeta\in\Pi_k}\left({\det\A^{jl}\over \det\A}\right)\left(k,\zeta\right)\right]~~(|\zeta|>0)
\\\\ \ds
~+~\sum_{k,l=1}^{n-N} \i^{v_j-u_l+1}\Big(g_l(t)-g_l(0)\Big)\left({\pi\i\over n-N}\right)\left({\det\A^{jl}\over\det\A}\right)(k,0)\qquad (\omega(0)=0)
\end{array}
\eeq
where $\Theta_k=\Big\{\vartheta\in[0,2\pi)\colon{(k-1)\pi\over n-N}<\vartheta<{k\pi\over n-N}\Big\}$ and $\Pi_k=\Big\{\vartheta\in[0,2\pi)\colon\vartheta={(k-1)\pi\over n-N}~\hbox{and}~\vartheta={k\pi\over n-N}\Big\}$,
for every $v_j\in\Big\{n-N,n-N+1,\ldots,n-1\Big\}$ and $L$ sufficiently large. }

\section{Novel Integral Representations}
\setcounter{equation}{0}
In this section, we develop the main framework of unified transform method, 
known as the {\it novel integral representations}.
It has been used for obtaining certain results in evolution P.D.E. See \cite{ASF1}-\cite{ASF2} and \cite{ASF3}.
An introduction can be found in \cite{ASF} and \cite{D.T.V}.

The PDE (\ref{PDE}) 
can be rewritten in the following divergence form:
\bel{divergence}
 \Big(e^{-\i x\xi+\omega(\xi)t}q(t,x)\Big)_t~-~\Bigg(e^{-\i x\xi+\omega(\xi)t}\sum_{m=0}^{n-1}
c_m(\xi)\partial_x^m q(t,x)\Bigg)_x~=~0,
\eeq
where $c_m(\xi)$ can be explicitly computed by
\bel{coe}
\sum_{m=0}^{n-1}c_m(\xi)\partial_x^m~=~\i{\omega(\xi)-\omega(\eta)\over \xi-\eta}\Bigg|_{\eta=-\i\partial_x}.
\eeq
From (\ref{coe}) and (\ref{omega_m}), we have
\bel{c_m}
c_m(\xi)~=~\i^{3m+1}\omega_{n-m-1}(\xi).
\eeq
Employing (\ref{divergence}) and using Green's theorem in the domain $\{0<x<\infty,~0<s<t\}$, we find the following {\it global relation}: 
\bel{Global}
e^{\omega(\xi)t}~\Hat{q}(t,\xi)~=~\Hat{q}_o(\xi)~-~\sum_{m=0}^{n-1}c_m(\xi)
\int_0^t e^{\omega(\xi)\tau}\partial_x^m q(\tau,0)d\tau
\eeq
which has an analytic continuation into the lower half plan for $\Im\xi\leq0$, 
where 
\bel{Fourier transform half-line}
\Hat{q}_o(\xi)~=~\int_0^\infty e^{-x\xi}q_o(x)dx\qquad\hbox{and}\qquad \Hat{q}(t,\xi)~=~\int_0^\infty e^{-x\xi} q(t,x)dx.
\eeq
We define the following $t$-transforms:
\bel{t transforms}
G_l(t,\xi)~=~\int_0^T e^{\omega(\xi)(\tau-t)}g_l(\tau)d\tau,~~~u_l\in U~~~\hbox{and}~~~ Q_j(t,\xi)~=~\int_0^T e^{\omega(\xi)(\tau-t)}\p_x^{v_j} q(0,\tau)d\tau,~~~ v_j\in V.
\eeq
By using the notations $\a_j=c_{v_j}$ and $\b_l=c_{u_l}$, equation (\ref{c_m}) implies   
\bel{a_j b_l}
\a_j(\xi)~=~\i^{3v_j+1}\omega_{n-v_j-1}(\xi),\qquad \b_l(\xi)~=~\i^{3u_l+1}\omega_{n-u_l-1}(\xi).
\eeq
Evaluating (\ref{Global}) at $t=T$ and multiplying the resulting equation by $e^{-\omega(\xi)t}$, we have
\bel{EQ}
\sum_{j=1}^N \a_{j}(\xi)Q_j(t,\xi)~=~-e^{\omega(\xi)(T-t)}~\Hat{q}(T,\xi)~+~e^{-\omega(\xi)t}~\Hat{q}_o(\xi)~-~\sum_{l=1}^{n-N}\b_{l}(\xi)G_l(t,\xi).
\eeq
By replacing $\xi$ in (\ref{EQ}) with $z_i^k(\xi)$, we obtain the following linear system of $N$ equations:
\bel{EQ_i}
\begin{array}{rl}\ds
\sum_{j=1}^N \a_{j}(z_i(\xi))Q_j(t,\xi)~=~-e^{\omega(\xi)(T-t)}~
\Hat{q}(T,z_i^k(\xi))~+~e^{-\omega(\xi)t}~\Hat{q}_o(z_i^k(\xi))~-~\sum_{l=1}^{n-N}\b_{l}(z_i^k(\xi))G_l(t,\xi),
\\ \ds
~~~~~~~~~~~~~~~~~~~~~~~~~~~~~~~~~~~~~~~~~~~~~~~~~~~ i=1,2,\ldots,N.
\end{array}
\eeq
The inverse of the matrix $\A$ defined in (\ref{A}) can be written as
\bel{A^-1}
\A^{-1}(k,\xi)~=~\left({\adj \A\over \det\A}\right)(k,\xi)
\eeq
where $\adj\A$ is the  {\it adjugate} of $\A$ whose $(i,j)$-th entry equals the determinant of the $(j,i)$-th principal minor of $\A$, denoted by $\A_{ji}$ multiplying $(-1)^{j+i}$. 
Recall from section 2 and (\ref{a_j b_l}). The matrix $\A^{jl}$ is constructed from $\A$ by replacing its $j$-th column by column vector $\Big(\b_l(z_i(\xi))\Big)_{1\leq i\leq N}$. 
By solving the system of the $N$-linear equations in (\ref{EQ_i}),  we have
\bel{Global relation}
\begin{array}{lr}\ds
Q_j(t,\xi)~=~-\i^{v_j-1}\sum_{i=1}^N(-1)^{j+i}\left({\det\A_{ji}\over \det\A}\right)(k,\xi)\left(e^{\omega(\xi)(T-t)}~\Hat{q}(T,z_i^k(\xi))-e^{-\omega(\xi)t}~\Hat{q}_o(z_i^k(\xi))\right)
\\\\ \ds ~~~~~~~~~~~~~~~~
~-~\i^{v_j-u_l}\sum_{l=1}^{n-N}\left({\det\A^{jl}\over \det\A}\right)(k,\xi) G_l(t,\xi),\qquad\qquad j=1,2,\ldots,N
\end{array}
\eeq
for every $k=1,2,\ldots,n-N$.
Suppose that we are given the canonical boundary values, i.e: $\U=\Big\{0,1,\ldots,n-N-1\Big\}$. Then,
\bel{det A}
\det\A(k,\xi)~=~a_n^N\prod_{1\leq i<j\leq N}\left(z_i^k(\xi)-z_j^k(\xi)\right)
\eeq
which is not vanishing for $\xi$ asymptotically large, as was discussed in \cite{FS}.

The expression in (\ref{det A}) can be observed as follows. By (\ref{A}), if $z_i=z_j$ for any pair of $i\neq j$, then $\det\A=0$. 
The determinant is a polynomial of $z_1, z_2,\ldots,z_N$ in an order of $N(N-1)/2$.  
Therefore, it must be devisable by the product in (\ref{det A}). On the other hand, the number of the factors in the product equals to $N(N-1)/2$ 
which is a combinatorial fact. 

Recall that $\A^{jl}$ is the matrix obtained from $\A$ after replacing its $j$-th column by $\omega_{n-u_l-1}(z_i)$ for every $i$-th row. 
If $z_i=z_j$ for any pair of $i\neq j$, then $\det\A^{jl}=0$. This implies that the determinant of $\A^{jl}$ is devisable by the product in (\ref{det A}) as well.
When $n\leq 5$, each $z_i=z_i(\xi),~i=1,2,\ldots,N$ can be computed explicitely by radical formulae, which is analytic for $\xi$ asymptotically large. 
From all above, the quotient $\det\A^{jl}/\det\A$ for every $j,l$ has possibly finitly many removable singularities, and is analytic for 
$\xi\in D^+_L$ provided that $L$ is sufficiently large. 

By multiplying equation (\ref{Global}) with $-\i\omega'(\xi)$  and integrating along the contour $\p D^+_L$ defined in (\ref{D_L union}),  
we have the following novel integral representation: 
\bel{NIT}
\begin{array}{lr}\ds
\int_{\p D^+_L}\left( \int_0^T- \i\omega'(\xi)e^{\omega(\xi)(\tau-t)}\p_x^{v_j}q(\tau,0)d\tau\right)d\xi~=~
\\\\ \ds~~~~~~~~~~~~~~~~~
\sum_{k=1}^{n-N}\int_{\p D^+_{L,k}}\sum_{i=1}^N\i^{v_j}(-1)^{j+i}\left({\det\A_{ji}\over \det\A}\right)(k,\xi)\omega'(\xi)\left(e^{\omega(\xi)(T-t)}~\Hat{q}(T,z_i(\xi))-e^{-\omega(\xi)t}~\Hat{q}_o(z_i(\xi))\right)d\xi
\\\\ \ds~~~~~~~~~~~~~~~
~+~\sum_{l=1}^{n-N} \i^{v_j-u_l+1}\left\{\sum_{k=1}^{n-N}\int_{\p D^+_{L,k}}\left({\det\A^{jl}\over \det\A}\right)(k,\xi)\omega'(\xi)\Bigg(\int_0^T e^{\omega(\xi)(\tau-t)}g_l(\tau)d\tau\Bigg)d\xi\right\}
\end{array}
\eeq
for every $j=1,2,\ldots,N$.

\section{Proof of Theorem Two}
\setcounter{equation}{0}
Our estimation will be carried out in the general setting where the given boundary conditions prescribed in (\ref{given conditions}) are not necessary in the canonical form. 
We assume that  
 \bel{assumption}
 \left({\det\A^{jl}\over \det\A}\right)(k,\xi)\qquad k,l=1,2,\ldots,n-N,\qquad j=1,2,\ldots,N
 \eeq
 have only finitely many removable singularities,
 and are analytic for $\xi$ asymptotically large.  If (\ref{given conditions}) is in the canonical form and $n\leq 5$, then this assumption can be removed.
The proof will be accomplished within several steps. 
 
{\bf 1.} Recall the definition of $\A_{ji}$, $\A^{jl}$ and (\ref{omega_m}) in section 3. 
Since $\omega(\xi)\approx a_n\xi^n$ for $\xi$ is asymptotically large, we shall have
\bel{Det est}
\left|\left({\det\A_{ji}\over\det\A}\right)(k,\xi)\right|~\lesssim~|\xi|^{v_j-n+1},\qquad\left| \left({\det\A^{jl}\over\det\A}\right)(k,\xi)\right|~\lesssim~|\xi|^{v_j-u_l}
\eeq
as $|\xi|\mt\infty$.
From (\ref{Poly})-(\ref{z}). We have $|z_i(\xi)|\sim|\xi|$  with $\Im{z_i^k(\xi)}<0$ as $|\xi|\mt\infty$ for $i=1,2,\ldots,N$ and $k=1,2,\ldots,n-N$. Consider the Fourier transform 
\bel{q hat}
\Hat{q}\left(T,z_i^k(\xi)\right)=\int_0^\infty e^{-\i xz_i^k(\xi)}q(T,x)dx
\eeq
where $q(T,x)$ has sufficiently decay as $x\mt\infty$. See \cite{ASF1}, \cite{ASF} and \cite{FS}.

Integration by parts with respect to $x$ inside (\ref{q hat}) gives
\bel{q Hat est}
\left|\Hat{q}\left(T, z_i^k(\xi)\right)\right|~\lesssim~|\xi|^{-1},\qquad \xi\in D^+_{L.k}.
\eeq
Recall the principal domain $\D$ in (\ref{principal domain}) and (\ref{prin domain +/-}). Define 
\bel{max bound}
\M(R)~=~\max_{|\xi|=R} ~\left\{\left|\Hat{q}\Big(T,z_i^k(\xi)\Big)\right|~\colon~\xi \in \D^+_k~,~i=1,2,\ldots,N~;~k=1,2,\ldots,n-N\right\}
\eeq
which tends to zero as $R\mt\infty$. 

{\bf 2.} By writing 
\[  a_n=\cos\varphi+\i\sin\varphi\qquad\hbox{and}\qquad \xi=R(\cos\vartheta+\i\sin\vartheta),\]     
we have $\cos(\varphi+n\vartheta)<0$ whereas $\varphi+n\vartheta\in\Big[{\pi\over 2},{3\pi\over 2}\Big]$ from (\ref{Degrees}), for $\xi\in\D$.

Consider the norm of
\bel{I_j}
\I_j~=~\sum_{k=1}^{n-N}\int_{\left\{\xi\in\D^+_k\colon |\xi|=R\right\}} 
\sum_{i=1}^N(-1)^{j+i}\left({\det\A_{ji}\over \det\A}\right)(k,\xi)\omega'(\xi)e^{\omega(\xi)(T-t)}~\Hat{q}\left(T,z_i^k(\xi)\right) d\xi.
\eeq
Let $|\xi|=R$. From the first inequality in (\ref{Det est}) we have 
\bel{bound set 1}
\begin{array}{lr}\ds
\M(R) \sum_{i=1}^N \left|\left({\det\A_{ji}\over \det\A}\right)  (k,\xi)\omega'(\xi)\right|
~\lesssim~\M(R) R^{v_j}\qquad\hbox{as}~~R\mt\infty.
\end{array}
\eeq
Recall that $\cos(\varphi+n\vartheta)<0$ is concave for $\vartheta\in\arg\D^+$ as defined in (\ref{Degrees}). For every $t<T$, 
 \bel{norm I_j}
\begin{array}{lr}\ds
\left|\I_j\right|~\lesssim~\sum_{k=1}^{n-N}\int_{\left\{\xi\in\D^+_k\colon |\xi|=R\right\}} 
\sum_{i=1}^N \left|\left({\det\A_{ji}\over \det\A}\right)  (k,\xi)\omega'(\xi)\Hat{q}\left(T,z_i^k(\xi)\right)\right|e^{\Re\omega(\xi)(T-t)} d\xi
\\\\ \ds~~~~~~
~\lesssim~\M(R) \int_{\arg\D^+} R^{v_j+1} e^{-R^n\cos(\varphi+n\vartheta)(T-t)}d\vartheta
\\\\ \ds~~~~~~
~\lesssim~\M(R) \int_{\arg\D^+} R^{v_j+1} e^{-R^n(T-t)\vartheta}d\vartheta
\\\\ \ds~~~~~~
~\lesssim~\M(R)\left({1\over R^{n-v_j-1}}\right)\left({1\over T-t}\right).
\end{array}
\eeq
Since $v_j$ is at most $n-1$ and $\M(R)$ turns to zero as $R$ approaching to infinity, 
we have $|\I|\mt0$ as $R\mt\infty$.     
Since $D^+$ in (\ref{Domain})-(\ref{D_L union}) coincides with $\D^+$ at infinity, by applying Cauchy's theorem on $D_L^+$ we have 
\bel{Est1}
\begin{array}{lr}\ds
\sum_{k=1}^{n-N}\int_{\p D^+_{L,k}} 
\sum_{i=1}^N(-1)^{j+i}\left({\det\A_{ji}\over \det\A}\right)(k,\xi)\omega'(\xi)e^{\omega(\xi)(T-t)}~\Hat{q}\left(T,z_i^k(\xi)\right) d\xi~=~0.
\end{array}
\eeq

{\bf 3.} Consider the left hand side of equation (\ref{NIT}). Since the integrant is analytic, we can deform back the contour from $\p D^+_L$ to $\p D^+$. Recall that $\p D=\{\xi\in\Cx\colon\Re{\omega(\xi)}=0\}$ and $D^+$ coincides with $\D^+$ at infinity which consist of $n-N$ sectors. 
By changing variables $\omega(\xi)=\i l$, we have $l\in\R$ and $-\i\omega'(\xi)d\xi=dl$. Therefore,  
\bel{p^v_j}
\begin{array}{lr}\ds
\int_{\p D^+}\left( \int_0^T- \i\omega'(\xi)e^{\omega(\xi)(\tau-t)}\p_x^{v_j}q(\tau,0)d\tau\right)d\xi~=~(n-N)\int_{-\infty}^\infty e^{il(\tau-t)}dl\left(\int_0^T \p_x^{v_j}q(\tau,0)d\tau\right)
\\\\
\qquad\qquad\qquad\qquad\qquad\qquad\qquad\qquad~~~~~\ds
~=~2\pi (n-N)\int_0^T\delta(\tau-t) \p_x^{v_j}q(\tau,0)d\tau
\\\\
\qquad\qquad\qquad\qquad\qquad\qquad\qquad\qquad~~~~~ \ds
~=~2\pi (n-N)\p_x^{v_j}q(t,0)
\end{array}
\eeq
 where $\delta$ is Dirac delta. 

{\bf 4.}  For the remaining terms in (\ref{NIT}), suppose $v_j<u_l$ and consider  
\bel{Integral separation}
\begin{array}{lr}\ds
\sum_{k=1}^{n-N}\int_{\p D^+_{L,k}}\left({\det\A^{jl}\over \det\A}\right)(k,\xi)
\omega'(\xi)\left(\int_0^T e^{\omega(\xi)(\tau-t)}g_l(\tau)d\tau\right)d\xi
\\\\ \ds
~=~\sum_{k=1}^{n-N}\int_{\p D^+_{L,k}}\left({\det\A^{jl}\over \det\A}\right)(k,\xi)
\omega'(\xi)\left(\int_0^t e^{\omega(\xi)(\tau-t)}g_l(\tau)d\tau+\int_t^T e^{\omega(\xi)(\tau-t)}g_l(\tau)d\tau\right)d\xi.
\end{array}
\eeq
Let
\bel{I_jl}
\I_{jl}~=~\sum_{k=1}^{n-N}\int_{\left\{\xi\in\D^+_k\colon |\xi|=R\right\}} \left({\det\A^{jl}\over \det\A}\right)(k,\xi)
\omega'(\xi) e^{\omega(\xi)(\tau-t)} d\xi.
\eeq
From the second inequality in (\ref{Det est}),  we have
\bel{bound set 2}
\left|\left({\det\A^{jl}\over \det\A}\right)(k,\xi)
\omega'(\xi)\right|~\lesssim~|\xi|^{v_j-u_l+(n-1)}\qquad \hbox{as}~~|\xi|\mt\infty.
\eeq
Recall that $\cos(\varphi+n\vartheta)<0$ is concave 
for $\vartheta\in\arg\D^+$. For every $\tau>t$, 
 \bel{norm I_l}
\begin{array}{lr}\ds
\left|\I_{jl}\right|~\lesssim~\sum_{k=1}^{n-N}\int_{\left\{\xi\in\D^+_k\colon |\xi|=R\right\}} 
\left|\left({\det\A^{jl}\over \det\A}\right) (k,\xi)\omega'(\xi)\right|e^{\Re\omega(\xi)(\tau-t)} d\xi
\\\\ \ds~~~~~~~
~\lesssim~ \int_{\arg\D^+} R^{v_j-u_l+n} e^{-R^n\cos(\varphi+n\vartheta)(\tau-t)}d\vartheta
\\\\ \ds~~~~~~~
~\lesssim~ \int_{\arg\D^+} R^{v_j-u_l+n} e^{-R^n(\tau-t)\vartheta}d\vartheta
\\\\ \ds~~~~~~~
~\lesssim~ R^{v_j-u_l}(\tau-t)^{-1}.
\end{array}
\eeq
By applying Cauchy's theorem on $D^+_L$, we have
 \bel{Est2} 
\sum_{k=1}^{n-N}  \int_{\p D^+_{L,k}}\left({\det\A^{jl}\over \det\A}\right)(k,\xi)
\omega'(\xi)\left(\int_t^T e^{\omega(\xi)(\tau-t)}g_l(\tau)d\tau\right)d\xi~=~0.
\eeq

{\bf 5.} From (\ref{principal domain})-(\ref{Degrees}), $\xi\in\D^+_k$ implies ${(k-1)\pi\over n-N}\leq\arg\xi\leq{k\pi\over n-N}$. Let
\bel{II_jl}
\II_{jl}~=~\sum_{k=1}^{n-N}  \int_{\left\{\xi\in\Cx^+\setminus\D^+_k\colon |\xi|=R, {(k-1)\pi\over n-N}\leq\arg\xi\leq{k\pi\over n-N}  \right\}}   
  \left({\det\A^{jl}\over \det\A}\right)(k,\xi)\omega'(\xi) e^{\omega(\xi)(\tau-t)} d\xi.
\eeq
We have $\cos(\varphi+n\vartheta)>0$ is convex for $\vartheta\in\arg\left(\Cx^+\setminus\D^+\right)$.
For every $\tau<t$, 
 \bel{norm II_l}
\begin{array}{lr}\ds
\left|\II_{jl}\right|~\lesssim~\sum_{k=1}^{n-N}\int_{\left\{\xi\in\Cx^+\setminus\D^+\colon |\xi|=R,~~{(k-1)\pi\over n-N}\leq\arg\xi\leq{k\pi\over n-N}\right\}} 
\left|\left({\det\A^{jl}\over \det\A}\right) (k,\xi)\omega'(\xi)\right|e^{\Re\omega(\xi)(\tau-t)} d\xi
\\\\ \ds~~~~~~~
~\lesssim~ \int_{\arg\left(\Cx^+\setminus\D^+\right)} R^{v_j-u_l+n} e^{-R^n\cos(\varphi+n\vartheta)(t-\tau)}d\vartheta
\\\\ \ds~~~~~~~
~\lesssim~ \int_{\arg\left(\Cx^+\setminus\D^+\right)} R^{v_j-u_l+n} e^{-R^n(t-\tau)\vartheta}d\vartheta
\\\\ \ds~~~~~~~
~\lesssim~R^{v_j-u_l}(t-\tau)^{-1}.
\end{array}
\eeq
By assumption, $\det\A^{jl}/\det\A$ has only finitely many removable singularities. By applying Cauchy's theorem on $\Cx^+\setminus D^+_L$, and Residue's theorem on $\left\{\xi\in D^+\colon |\xi|\leq L\right\}$,
we have
\bel{Est2'}
\begin{array}{lr}\ds
  \sum_{k=1}^{n-N}  \int_{\p D^+_{L_k}}\left({\det\A^{jl}\over \det\A}\right)(k,\xi)\omega'(\xi)\left(\int_0^t e^{\omega(\xi)(\tau-t)}g_l(\tau)d\tau\right)d\xi 
   \\\\ \ds~~~~~~~
 ~=~ \sum_{k=1}^{n-N} \int_{0}^{\infty}\left[ \left({\det\A^{jl}\over \det\A}\right)\left(k,\lambda^k\xi\right)- \left({\det\A^{jl}\over \det\A}\right)\left(k,\lambda^{k-1}\xi\right) \right]  \omega'(\xi)\left(\int_0^t e^{\omega(\xi)(\tau-t)}g_l(\tau)d\tau\right)d\xi
\end{array}
   \eeq
   where $\lambda$ is defined in (\ref{lambda}).

 {\bf 6.} Suppose $v_j>u_l$ and consider
\bel{t-separation}
\begin{array}{lr}\ds
\sum_{k=1}^{n-N} \int_{\p D^+_{L,k}}\left({\det\A^{jl}\over \det\A}\right)(k,\xi)\omega'(\xi)\Bigg(\int_0^T e^{\omega(\xi)(\tau-t)}g_l(\tau)d\tau\Bigg)d\xi
\\\\ \ds
~=~\sum_{k=1}^{n-N}\int_{\p D^+_{L,k}}\left({\det\A^{jl}\over \det\A}\right)(k,\xi)\omega'(\xi)\Bigg(\int_t^T e^{\omega(\xi)(\tau-t)}g_l(\tau)d\tau+\int_0^t e^{\omega(\xi)(\tau-t)}g_l(\tau)d\tau\Bigg)d\xi
\end{array}
\eeq
Integrating by parts with respect to $\tau$ yields 
\bel{int by parts}
\begin{array}{lr}\ds
\sum_{k=1}^{n-N}\int_{\p D^+_{L,k}}\left({\det\A^{jl}\over \det\A}\right)(k,\xi){\omega'(\xi)\over \omega(\xi)}\left( e^{\omega(\xi)(T-t)}g_l(T)-\int_t^T e^{\omega(\xi)(\tau-t)}\dot{g}_l(\tau)d\tau\right)d\xi
\\\\ \ds
~-~\sum_{k=1}^{n-N}\int_{\p D^+_{L,k}}\left({\det\A^{jl}\over \det\A}\right)(k,\xi){\omega'(\xi)\over \omega(\xi)}\left(e^{-\omega(\xi)t}g_l(0)+\int_0^t e^{\omega(\xi)(\tau-t)}\dot{g}_l(\tau)d\tau\right)d\xi.
\end{array}
\eeq
Notice that $\omega'(\xi)/\omega (\xi)$ has possibly simple poles which can occur only on the boundary $\p D$ whereas $\Re\omega(\xi)=0$.

From the second inequality in (\ref{Det est}), we have 
\bel{bound est 3}
\begin{array}{lr}\ds
\left|\left({\det\A^{ji}\over \det\A}\right)  (k,\xi){\omega'(\xi)\over\omega(\xi)}\right|
~\lesssim~|\xi|^{v_j-u_l-1}\qquad\hbox{as}~~|\xi|\mt\infty.
\end{array}
\eeq
Let
\bel{III_jl}
\III_{jl}~=~\sum_{k=1}^{n-N} \int_{\left\{\xi\in\D^+\colon |\xi|=R\right\}} \left({\det\A^{jl}\over \det\A}\right)(k,\xi){\omega'(\xi)\over\omega(\xi)} e^{\omega(\xi)(\tau-t)} d\xi.
\eeq
Recall that $\cos(\varphi+n\vartheta)<0$ is concave for $\vartheta\in\arg\D^+$. 
For every $\tau>t$, 
 \bel{norm III_l}
\begin{array}{lr}\ds
\left|\III_{jl}\right|~\lesssim~\sum_{k=1}^{n-N}\int_{\left\{\xi\in\D^+\colon |\xi|=R\right\}} 
\left|\left({\det\A^{jl}\over \det\A}\right) (k,\xi){\omega'(\xi)\over\omega(\xi)}\right|e^{\Re\omega(\xi)(\tau-t)} d\xi
\\\\ \ds~~~~~~~~~
~\lesssim~ \int_{\arg\D^+} R^{v_j-u_l} e^{R^n\cos(\varphi+n\vartheta)(\tau-t)}d\vartheta
\\\\ \ds~~~~~~~~~
~\lesssim~ \int_{\arg\D^+} R^{v_j-u_l} e^{-R^n(\tau-t)\vartheta}d\vartheta
\\\\ \ds~~~~~~~~~
~\lesssim~R^{v_j-u_l-n}(\tau-t)^{-1}.
\end{array}
\eeq
The result holds for $\tau-t$ replaced by $T-t$. 

By applying Cauchy's theorem on $D^+_L$, we have
\bel{Est3}
\sum_{k=1}^{n-N} \int_{\p D^+_{L,k}}\left({\det\A^{jl}\over \det\A}\right)(k,\xi){\omega'(\xi)\over \omega(\xi)}\left( e^{\omega(\xi)(T-t)}g_l(T)- \int_t^T e^{\omega(\xi)(\tau-t)}\dot{g}_l(\tau)d\tau\right)d\xi~=~0.
\eeq

{\bf 7.} Let 
\bel{IV_jl}
\IV_{jl}~=~ \sum_{k=1}^{n-N} \int_{\left\{\xi\in\Cx^+\setminus\D^+_k\colon |\xi|=R,{(k-1)\pi\over n-N}\leq\arg\xi\leq{k\pi\over n-N}\right\}}   
   \left({\det\A^{jl}\over \det\A}\right)(k,\xi){\omega'(\xi)\over\omega(\xi)} e^{\omega(\xi)(\tau-t)} d\xi.
\eeq
Recall that $\cos(\varphi+n\vartheta)>0$ is convex for $\vartheta\in\left(\arg\Cx^+\setminus\D^+\right)$. 
For every $\tau<t$,
 \bel{norm IV_jl}
\begin{array}{lr}\ds
\left|\IV_{jl}\right|~\lesssim~\sum_{k=1}^{n-N}\int_{\left\{\xi\in\Cx^+\setminus\D^+\colon |\xi|=R,{(k-1)\pi\over n-N}\leq\arg\xi\leq{k\pi\over n-N}\right\}} 
\left|\left({\det\A^{jl}\over \det\A}\right) (k,\xi){\omega'(\xi)\over\omega(\xi)}\right|e^{\Re\omega(\xi)(\tau-t)} d\xi
\\\\ \ds~~~~~~~~
~\lesssim~ \int_{\arg\left(\Cx^+\setminus\D^+\right)} R^{v_j-u_l} e^{-R^n\cos(\varphi+n\vartheta)(t-\tau)}d\vartheta
\\\\ \ds~~~~~~~~
~\lesssim~ \int_{\arg\left(\Cx^+\setminus\D^+\right)} R^{v_j-u_l} e^{-R^n(t-\tau)\vartheta}d\vartheta
\\\\ \ds~~~~~~~~
~\lesssim~R^{v_j-u_l-n}(t-\tau)^{-1}.
\end{array}
\eeq
By applying Cauchy's theorem on $\Cx^+\setminus D^+_L$, and  Residue's theorem on $\left\{\xi\in D^+\colon |\xi|\leq L\right\}$,
we have
\bel{Est3'}
\begin{array}{lr}\ds
   \sum_{k=1}^{n-N}\int_{\p D^+_{L,k}}\left({\det\A^{jl}\over \det\A}\right)(k,\xi){\omega'(\xi)\over \omega(\xi)}\left( \int_0^t e^{\omega(\xi)(\tau-t)}\dot{g}_l(\tau)d\tau\right)d\xi
   \\\\ \ds ~~~~~~~
~=~\sum_{k=1}^{n-N}\hbox{\bf p.v}\int_{0}^{\infty}\left[\left({\det\A^{jl}\over \det\A}\right)\left(k,\lambda^k\xi\right)-\left({\det\A^{jl}\over \det\A}\right)\left(k,\lambda^{k-1}\xi\right)\right]{\omega'(\xi)\over \omega(\xi)}\left(\int_0^t e^{\omega(\xi)(\tau-t)}\dot{g}_l(\tau)d\tau\right)d\xi
\\\\ \ds ~~~~~~~
~+~\sum_{k=1}^{n-N} \i^{v_j-u_l+1}\Big(g_l(t)-g_l(0)\Big)\left[ 2\pi\i\sum_{\omega(\zeta)=0,\arg\zeta\in\Theta_k} \left({\det\A^{jl}\over \det\A}\right)\left(k,\zeta\right)+\pi\i\sum_{\omega(\zeta)=0,\arg\zeta\in\Pi_k}\left({\det\A^{jl}\over \det\A}\right)\left(k,\zeta\right)\right]~~(|\zeta|>0)
\\\\ \ds~~~~~~~
~+~\sum_{k=1}^{n-N} \i^{v_j-u_l+1}\Big(g_l(t)-g_l(0)\Big)\left({\pi\i\over n-N}\right)\left({\det\A^{jl}\over\det\A}\right)(k,0)\qquad (\omega(0)=0)
\end{array}
\eeq
where
\bel{Theta_k}
\Theta_k~=~\left\{\vartheta\in[0,2\pi)\colon {(k-1)\pi\over n-N}<\vartheta<{k\pi\over n-N}\right\}
\eeq
and
\bel{Pi_k}
\Pi_k~=~\left\{\vartheta\in[0,2\pi)\colon \vartheta={(k-1)\pi\over n-N}~~\hbox{and}~~\vartheta={k\pi\over n-N}\right\}
\eeq
for every $k=1,2,\ldots,n-N$.
\begin{remark}
In (\ref{Est2'}) and (\ref{Est3'}), we deform the contour $\p D^+_{L,k}$ to
\bel{Deform}
\left\{\arg\xi={k\pi\over n-N}\right\}\cup\left\{|\xi|=L,{(k-1)\pi\over n-N}<\arg\xi<{k\pi\over n-N}\right\}\cup\left\{\arg\xi={(k-1)\pi\over n-N}\right\}
\eeq
for every $k=1,2,\ldots,n-N$ with the same orientation.
\end{remark}

Lastly, recall that $g_l(0)=q_o(0)$ for every $u_l\in U$. 
By combining all above estimates, we have
\bel{Est}
\begin{array}{lr}\ds
2\pi(n-N)\p_x^{v_j}q(t,0)~=~-\sum_{k=1}^{n-N}\i^{v_j}\int_{\p D^+_{L,k}} \sum_{i=1}^N (-1)^{i+j}\left({\det \A_{ji}\over \det\A}\right)(k,\xi)\omega'(\xi)\Hat{q}_o\left(z_i^k(\xi)\right)e^{-\omega(\xi)t}d\xi
\\\\ \ds
~-~\sum_{k=1}^{n-N}\int_{\p D^+_{L,k}} \sum_{u_l<v_j} \i^{v_j-u_l+1}\left({\det\A^{jl}\over \det\A}\right)(k,\xi){\omega'(\xi)\over \omega(\xi)}q_o(0) 
e^{-\omega(\xi)t}d\xi
\\\\ \ds
~+~\sum_{u_l>v_j}\i^{v_j-u_l+1}\sum_{k=1}^{n-N} \int_{0}^{\infty}\left[ \left({\det\A^{jl}\over \det\A}\right)\left(k,\lambda^k\xi\right)- \left({\det\A^{jl}\over \det\A}\right)\left(k,\lambda^{k-1}\xi\right) \right]  \omega'(\xi)\left(\int_0^t e^{\omega(\xi)(\tau-t)}g_l(\tau)d\tau\right)d\xi\\\\ \ds
~+~ \sum_{u_l<v_j} \i^{v_j-u_l+1} \sum_{k=1}^{n-N}\hbox{\bf p.v}\int_{0}^{\infty}\left[\left({\det\A^{jl}\over \det\A}\right)\left(k,\lambda^k\xi\right)-\left({\det\A^{jl}\over \det\A}\right)\left(k,\lambda^{k-1}\xi\right)\right]{\omega'(\xi)\over \omega(\xi)}\left(\int_0^t e^{\omega(\xi)(\tau-t)}\dot{g}_l(\tau)d\tau\right)d\xi
\\\\ \ds
~+~\sum_{u_l<v_j}\sum_{k=1}^{n-N} \i^{v_j-u_l+1}\Big(g_l(t)-g_l(0)\Big)\left[ 2\pi\i\sum_{\omega(\zeta)=0,\arg\zeta\in\Theta_k} \left({\det\A^{jl}\over \det\A}\right)\left(k,\zeta\right)+\pi\i\sum_{\omega(\zeta)=0,\arg\zeta\in\Pi_k}\left({\det\A^{jl}\over \det\A}\right)\left(k,\zeta\right)\right]~~(|\zeta|>0)
\\\\ \ds
~+~\sum_{u_l<v_j}\sum_{k=1}^{n-N} \i^{v_j-u_l+1}\Big(g_l(t)-g_l(0)\Big)\left({\pi\i\over n-N}\right)\left({\det\A^{jl}\over\det\A}\right)(k,0)\qquad (\omega(0)=0).
\end{array}
\eeq
\endproof

\section{Proof of Theorem One}
\setcounter{equation}{0}
Let $\omega(\xi)=a_n\xi^n$. Domain $D$ in (\ref{Domain}) is replaced by the principal domain $\D$ in (\ref{principal domain})-(\ref{Degrees}), which is an union of $n$ sectors, as oriented in  (\ref{prin domain +/-}).
The following quotients of determinants can be computed explicitly : 
 \bel{Det explicit1}
\left({\det\A_{ji}\over\det\A}\right)(k,\xi)~=~a_n^{-1}\left({\det\V_{ji}\over\det\V}\right)(\rho^{n-N+1-k})\xi^{v_j-n+1},\qquad\qquad\xi\in\D^+_k,
\eeq
\bel{Det explicit2}
 \left({\det\A^{jl}\over\det\A}\right)(k,\xi)~=~\left({\det\V^{jl}\over\det\V}\right)(\rho^{n-N+1-k})\xi^{v_j-u_l},\qquad\qquad\xi\in\D^+_k.
 \eeq
By permutations, we can assume that $z_j\in\D^-_j$. By definition in (\ref{z}), we have  
\bel{rota}
z_j^k(\xi)~=~\rho^{j+n-N-k} \xi~\in~\D^-_{j}\qquad\hbox{whenever}\qquad \xi\in\D^+_k
\eeq
for every $k=1,2,\ldots,n-N$, and every $j=1,2,\ldots,N$.

From (\ref{Est3'}) in the previous section, the contour $D^+_L$ is replaced by $\D^+$, 
together with the explicit expressions in (\ref{Det explicit1})-(\ref{Det explicit2}) and (\ref{rota}).  

In particular, we have
\bel{Est explicit 1}
\begin{array}{lr}\ds
\sum_{k=1}^{n-N} \i^{v_j}\int_{\p \D^+_k} \sum_{i=1}^N (-1)^{i+j}\left({\det \A_{ji}\over \det\A}\right)(k,\xi)\omega'(\xi)\Hat{q}_o\left(z_i^k(\xi)\right)e^{-\omega(\xi)t}d\xi
\\\\ \ds
~=~\sum_{k=1}^{n-N}  \int_{\p \D^+_k}\sum_{i=1}^N(-1)^{i+j}\left({\det \V_{ji}\over \det\V}\right)(\rho^{n-N+1-k})(\i\xi)^{v_j}n\Hat{q}_o\left(\rho^{n-N+i-k}\xi\right)e^{-\omega(\xi)t}d\xi
\\\\ \ds
\sum_{k=1}^{n-N}\int_{\p\D^+_k} \sum_{u_l<v_j} \i^{-u_l+1}\left({\det\A^{jl}\over \det\A}\right)(k,\xi){\omega'(\xi)\over \omega(\xi)}q_o(0)e^{-\omega(\xi)t}d\xi
\\\\ \ds
~=~
-\sum_{k=1}^{n-N}  \int_{\p \D^+}\sum_{u_l<v_j} \left({\det\V^{jl}\over\det\V}\right)(\rho^{n-N+1-k})(\i\xi)^{v_j-u_l-1}nq_o(0)e^{-\omega(\xi)t}d\xi.
\end{array}
\eeq
On the other hand,  from (\ref{rota}) we have $z_j^k(\rho\xi)=z_j^{k-1}(\xi)$ for every $k=1,2,\ldots, n-N$ which implies that integrations over every $\p\D^+_k$ are equivalent. Recall the rays $\mathcal{R}^1_k$ and $\mathcal{R}^2_k$ defined below (\ref{prin domain boundary}).

For $v_j<u_l$, we have
 \bel{Est explicit 2}
 \begin{array}{lr}\ds
 \int_{0}^\infty\left[\left({\det\A^{jl}\over \det\A}\right)(n-N,\xi)  -\left({\det\A^{jl}\over \det\A}\right)(1,\xi) \right] \omega'(\xi) 
 \left(\int_0^t e^{\omega(\xi)(\tau-t)}g_l(\tau)d\tau\right)d\xi 
\\\\ \ds \qquad\qquad
~=~\left({\det\V^{jl}\over\det\V}\right)(\rho)\int_0^t \left\{\int_{\mathcal{R}^2_{n-N}} \xi^{v_j-u_l}\omega'(\xi)e^{\omega(\xi)(\tau-t)}d\xi\right\} g_l(\tau)d\tau
 \\\\\ds \qquad\qquad
 ~+~\left({\det\V^{jl}\over\det\V}\right)(\rho^{n-N})\int_0^t \left\{\int_{\mathcal{R}^1_{1}} \xi^{v_j-u_l}\omega'(\xi)e^{\omega(\xi)(\tau-t)}d\xi\right\} g_l(\tau)d\tau 
 \end{array}
 \eeq
For $v_j>u_l$, we have
\bel{Est explicit 3}
 \begin{array}{lr}\ds
 \int_{0}^\infty\left[\left({\det\A^{jl}\over \det\A}\right)(n-N,\xi)  -\left({\det\A^{jl}\over \det\A}\right)(1,\xi) \right] {\omega'(\xi)\over\omega(\xi)}\left(\int_0^t e^{\omega(\xi)(\tau-t)}\dot{g}_l(\tau)d\tau\right)d\xi 
\\\\ \ds \qquad\qquad
~=~\left({\det\V^{jl}\over\det\V}\right)(\rho)\int_0^t \left\{\int_{\mathcal{R}^2_{n-N}} \xi^{v_j-u_l}{\omega'(\xi)\over\omega(\xi)}e^{\omega(\xi)(\tau-t)}d\xi\right\} \dot{g}_l(\tau)d\tau
 \\\\\ds \qquad\qquad
 ~+~\left({\det\V^{jl}\over\det\V}\right)(\rho^{n-N})
\int_0^t \left\{\int_{\mathcal{R}^1_{1}} \xi^{v_j-u_l}{\omega'(\xi)\over\omega(\xi)}e^{\omega(\xi)(\tau-t)}d\xi\right\} \dot{g}_l(\tau)d\tau. 
\end{array}
 \eeq
Let $\eta=\omega(\xi)$ and write $\eta=Re^{\i\theta}$ in polar coordinates. From (\ref{vartheta Est}), 
\bel{variable change1}
\omega(\xi)~=~\eta~=~R e^{\i{\pi\over 2}}\qquad\hbox{for}~~\xi\in\mathcal{R}^1_1
\eeq
and
\bel{variable change2}
\omega(\xi)~=~\eta~=~Re^{\i{3\pi\over 2}}\qquad\hbox{for}~~\xi\in\mathcal{R}^2_{n-N}.
\eeq
From the estimation in the previous section,  we can deform the contours by rotating $\mathcal{R}^1_1$ by ${\pi\over 2n}$ clockwisely, and by rotating $\mathcal{R}^2_{n-N}$ by ${\pi\over 2n}$ counterclockwisely. 

Let
 \bel{tilde vartheta}
 \Tilde{\vartheta}_1~=~\vartheta_1-{\pi\over 2n},\qquad\qquad\Tilde{\vartheta}_2~=~\vartheta_2+{\pi\over 2n}.
 \eeq
 
 Equation  (\ref{Est explicit 2}) becomes
 \bel{Est explicit 2'}
  \begin{array}{lr}\ds
 \left({\det\V^{jl}\over\det\V}\right)(\rho)e^{\i(v_j-u_l)\Tilde{\vartheta}_2}\int_0^t \left(\int_{\infty}^0  \eta^{v_j-u_l\over n}e^{\eta(\tau-t)}d\eta\right) g_l(\tau)d\tau
 \\\\\ds \qquad
 ~+~\left({\det\V^{jl}\over\det\V}\right)(\rho^{n-N})e^{\i(v_j-u_l)\Tilde{\vartheta}_1}\int_0^t \left(\int_{0}^\infty \eta^{v_j-u_l\over n}e^{\eta(\tau-t)}d\eta\right) g_l(\tau)d\tau .
 \end{array}
 \eeq
 
Similarly, equation (\ref{Est explicit 3}) becomes
\bel{Est explicit 3'}
\begin{array}{lr}\ds
\left({\det\V^{jl}\over\det\V}\right)(\rho)e^{\i(v_j-u_l)\Tilde{\vartheta}_2}\left\{\int_0^t \left(\int_{\infty}^0 \eta^{v_j-u_l-n\over n}e^{\eta(\tau-t)}d\eta\right) \dot{g}_l(\tau)d\tau\right\}
 \\\\\ds ~
 ~+~\left({\det\V^{jl}\over\det\V}\right)(\rho^{n-N})e^{\i(v_j-u_l)\Tilde{\vartheta}_1}\left\{\int_0^t \left(\int_{0}^\infty \eta^{v_j-u_l-n\over n}e^{\eta(\tau-t)}d\eta\right) \dot{g}_l(\tau)d\tau\right\}. 
\end{array}
\eeq

By letting $-s=\eta(\tau-t)$, we find
\bel{Gamma1}
\begin{array}{lr}\ds
\int_0^t\left(\int_0^\infty \eta^{v_j-u_l\over n}e^{\eta(\tau-t)}d\eta\right) g_l(\tau)d\tau~=~\int_0^\infty s^{v_j-u_l\over n}e^{-s}ds \cdot\int_0^t {g_l(\tau)d\tau\over (t-\tau)^{v_j-u_l+n\over n}}
\\\\ \ds ~~~~~~~~~~~~~~~~~~~~~~~~~~~~~~~~~~~~~~~~~~~~~~~~~~~~~~
~=~\Gamma\left({v_j-u_l+n\over n}\right)\int_0^t {g_l(\tau)d\tau\over (t-\tau)^{v_j-u_l+n\over n}}
\end{array}
\eeq

and simultaneously
\bel{Gamma2}
\begin{array}{lr}\ds
\int_0^t\left(\int_0^\infty \eta^{v_j-u_l+n\over n}e^{\eta(\tau-t)}d\eta\right) \dot{g}_l(\tau)d\tau~=~\int_0^\infty s^{v_j-u_l-n\over n}e^{-s}ds \cdot\int_0^t {\dot{g}_l(\tau)d\tau\over (t-\tau)^{v_j-u_l\over n}}
\\\\ \ds ~~~~~~~~~~~~~~~~~~~~~~~~~~~~~~~~~~~~~~~~~~~~~~~~~~~~~~~~
~=~\Gamma\left({v_j-u_l \over n}\right)\int_0^t {\dot{g}_l(\tau)d\tau\over (t-\tau)^{v_j-u_l\over n}}.
\end{array}
\eeq

By bringing the above estimates into (\ref{Est}), we have
 \bel{Est Explicit}
\begin{array}{lr}\ds
 2\pi\p^{v_j}_x q(t,0)~=~-{1\over n-N}\sum_{k=1}^{n-N}\int_{\p \D^+_k}\sum_{i=1}^N(-1)^{i+j}\left({\det \V_{ji}\over \det\V}\right)(\rho^{n-N+1-k})n\Hat{q}_o(\rho^{n-N+i-k}\xi)(\i\xi)^{v_j}
 \\\\ \ds
 ~+~{1\over n-N}\sum_{k=1}^{n-N} \int_{\p \D^+_k}\sum_{u_l<v_j}\left({\det\V^{jl}\over\det\V}\right)(\rho^{n-N+1-k})nq_o(0)(\i\xi)^{v_j-u_l-1}e^{-\omega(\xi)t}d\xi
 \\\\ \ds 
 ~+~ \sum_{v_j<u_l}\i^{v_j-u_l+1} \Lambda_{jl}\Gamma\left({v_j-u_l+n\over n}\right)\int_0^t {g_l(\tau)d\tau\over (t-\tau)^{v_j-u_l+n\over n}}
 ~-~\sum_{u_l<v_j}\i^{v_j-u_l+1}\Lambda_{jl}\Gamma\left({v_j-u_l\over n}\right) \int_0^t {\dot{g}_l(\tau)d\tau\over (t-\tau)^{v_j-u_l\over n}} 
 \end{array}
 \eeq
 for every $v_j\in V$, where  
\bel{Lambda_jl}
\Lambda_{jl}~=~\left({\det\V^{jl}\over\det\V}\right)(\rho^{n-N})\exp\left(\i(v_j-u_l)\left(\vartheta_1^1-{\pi\over 2n}\right)\right)-\left({\det\V^{jl}\over\det\V}\right)(\rho)\exp\left(\i(v_j-u_l)\left(\vartheta_{n-N}^2+{\pi\over 2n}\right)\right)
\eeq
with $\vartheta_1$ and $\vartheta_2$ defined in (\ref{vartheta Est}).
\endproof

\section{Examples}
\setcounter{equation}{0}
The last section is devoted to a number of examples. For brevity, we assume the zero initial condition: $q_o(x)\equiv0$ for $x\in(0,\infty)$. 

 {\bf Example 1.} Our first example is the heat equation: 
 \bel{heat}
q_t(t,x)-q_{xx}(t,x)~=~0,\qquad 0<x<\infty,\qquad 0<t<T.
\eeq
Regarding to (\ref{heat}), we have $a_n=1$,  $n=2$ and $N=1$. The principal domain $\D$ is shown as below where the rotation $\rho=e^{\i\pi}$ maps $\D^+$ to $\D^-$.
 \begin{figure}[h]
\centering
\includegraphics[scale=0.40]{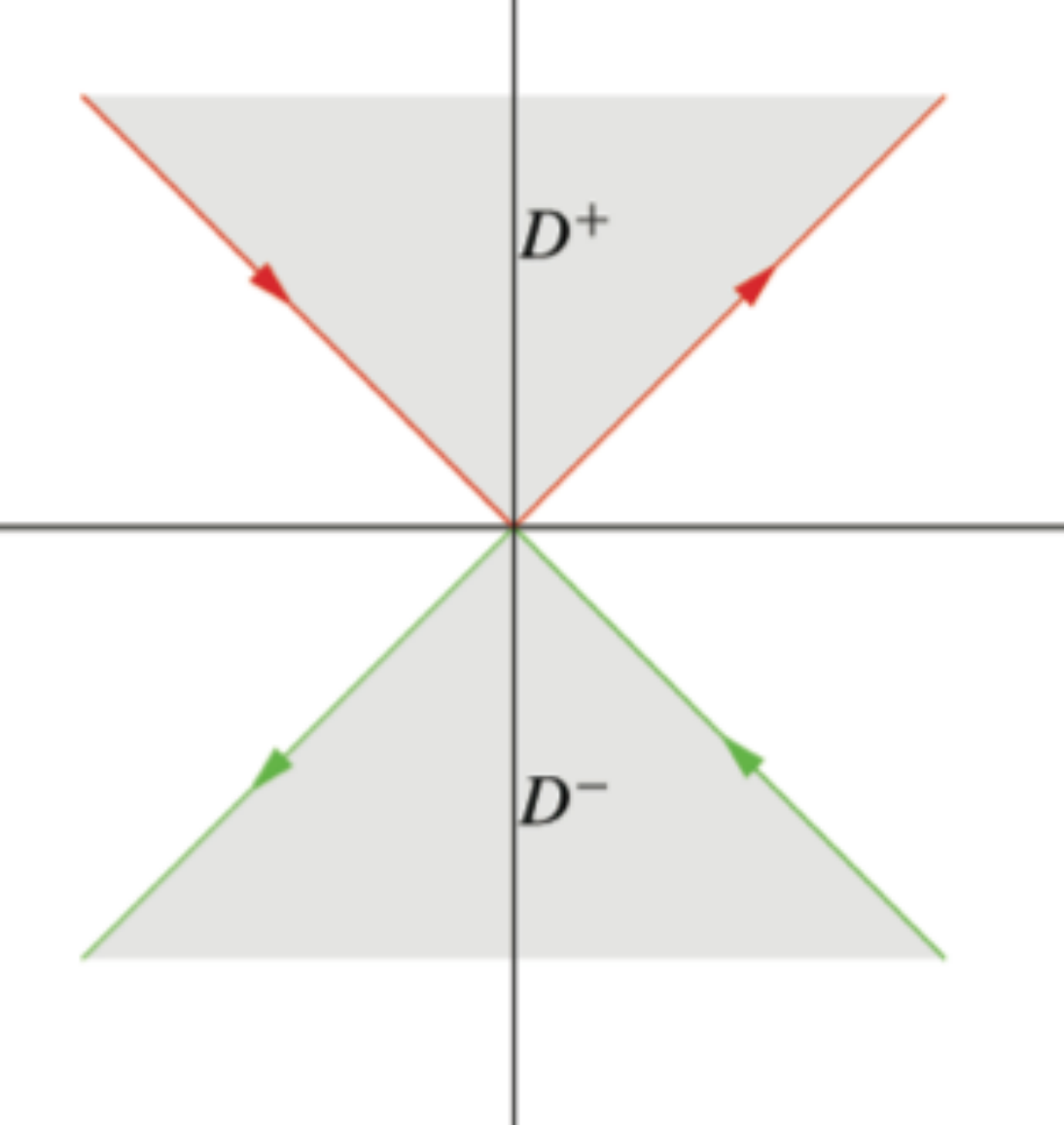}
\caption{$\p\D^+$ is represented as the red contour .}
\end{figure}

From (\ref{vartheta Est}), by taking $n=2$ and $\varphi=0$, we find 
\bel{vartheta(1)}
\vartheta_1^1~=~{\pi\over 4},\qquad \vartheta_1^2~=~{3\pi\over 4}.
\eeq
Suppose that the Dirichlet condition $q(t,0)=g(t)$ is given. We aim to find the Neumann boundary value $q_x(t,0)$ in terms of $g(t)$. In this case, we have $v_1=1$, $u_1=0$. From (\ref{V}), the matrix $\V(\rho)$ is a scalar equal to $1$ and $\V^{11}(\rho)=-1$.  

Bringing all these values into formula (\ref{Lambda_jl}), we find 
\bel{Lambda_11(1)}
\begin{array}{lr}\ds
\Lambda_{11}~=~
\left({\det\V^{11}\over\det\V}\right)(\rho)\exp\left(\i(v_1-u_1)\left(\vartheta_1^1-{\pi\over 4}\right)\right)-\left({\det\V^{11}\over\det\V}\right)(\rho)\exp\left(\i(v_1-u_1)\left(\vartheta_1^2+{\pi\over 4}\right)\right)
\\\\ \ds~~~~~~
~=~-1\cdot 1~-~(-1)(-1)~=~-2.
\end{array}
\eeq
By Theorem One, we have
\bel{result1}
q_x(t,0)~=~{-1\over\pi}\Gamma\left({1\over 2}\right)\int_0^t {\dot{g}(\tau)d\tau\over \sqrt{t-\tau}}.
\eeq
On the other hand, let Neumann condition $q_x(t,0)=g(t)$ be given. We aim to find $q(t,0)$ in terms of $g(t)$. 
We then have $v_1=1$, $u_1=0$. From (\ref{V}), the matrix $\V(\rho)$ is a scalar equal to $-1$ and $\V^{11}(\rho)=1$. 
Direct computation shows that $\Lambda_{11}=-2$ as in (\ref{Lambda_11(1)}). 

By Theorem One, we have 
\bel{result1'}
q(t,0)~=~{-1\over\pi}\Gamma\left({1\over 2}\right)\int_0^t {g(\tau)d\tau\over \sqrt{t-\tau}}.
\eeq

{\bf Example 2.} Consider a linear evolution equation with a third order derivative:
\bel{Third1}
q_t(t,x)+q_{xxx}(t,x)~=~0,\qquad0<x<\infty,\qquad0<t<T.
\eeq
Regarding to (\ref{Third1}), we have $a_n=-\i$, $n=3$ and $N=2$. The principal domain $\D$ is shown as below, where the rotation $\rho=e^{2\pi\i\over 3}$ maps $\D^+$ to $\D^-_1$ and $\rho^2$ maps $\D^+$ to $\D^-_2$.
 \begin{figure}[h]
\centering
\includegraphics[scale=0.40]{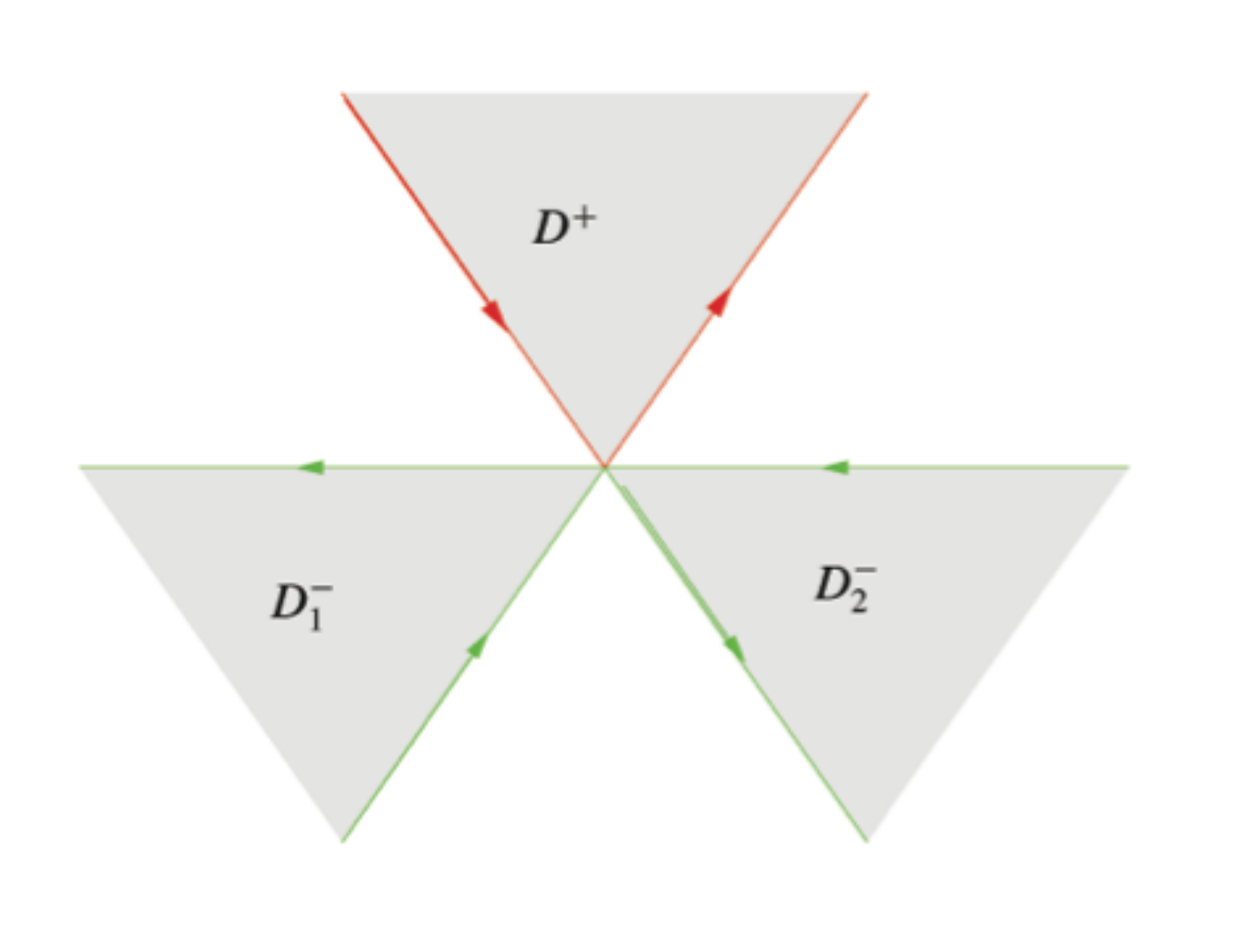}
\caption{$\p\D^+$ is represented as the red contour .}
\end{figure}

By taking $n=3$ and $\varphi={3\pi\over 2}$ in (\ref{vartheta Est}), we find 
\bel{vartheta(2)}
\vartheta_1^1~=~{\pi\over 3}\qquad\hbox{and}\qquad \vartheta_1^2~=~{2\pi\over 3}.
\eeq
Let $q(t,0)=g(t)$ be given. We aim to find $q_x(t,0)$ and $q_{xx}(t,0)$ in terms of $g(t)$. 
In this case, we have $u_1=0$, $v_1=1$ and $v_2=2$. Therefore, 
\bel{V example2}
\V(\rho)~=~\left[\begin{array}{lr}\ds
\rho\qquad1
\\\\ \ds
\rho^2\qquad 1
\end{array}\right],\qquad 
\V^{11}(\rho)~=~\left[\begin{array}{lr}\ds
\rho^2\qquad1
\\\\ \ds
\rho\qquad 1
\end{array}\right],\qquad 
\V^{21}(\rho)~=~\left[\begin{array}{lr}\ds
\rho\qquad\rho^2
\\\\ \ds
\rho^2\qquad \rho
\end{array}\right].
\eeq
From computations, we have
\bel{quotients detV example2}
{\det\V^{11}\over\det\V}(\rho)~=~ {\det\V^{21}\over\det\V}(\rho)~=~-1.
\eeq
 By (\ref{Lambda_jl}), we have 
\bel{Lambda_11(2)}
\begin{array}{lr}\ds
\Lambda_{11}~=~
\left({\det\V^{11}\over\det\V}\right)(\rho)\exp\left(\i(v_1-u_1)\left(\vartheta_1^1-{\pi\over 6}\right)\right)-\left({\det\V^{11}\over\det\V}\right)(\rho)\exp\left(\i(v_1-u_1)\left(\vartheta_1^2+{\pi\over 6}\right)\right)
\\\\ \ds~~~~~~
~=~(-1)\exp\left(\i{\pi\over 6}\right)~-~(-1)\exp\left(\i{5\pi\over 6}\right)~=~-\sqrt{3},
\end{array}
\eeq
and
\bel{Lambda_21(2)}
\begin{array}{lr}\ds
\Lambda_{21}~=~
\left({\det\V^{21}\over\det\V}\right)(\rho)\exp\left(\i(v_2-u_1)\left(\vartheta_1^1-{\pi\over 6}\right)\right)-\left({\det\V^{21}\over\det\V}\right)(\rho)\exp\left(\i(v_2-u_1)\left(\vartheta_1^2+{\pi\over 6}\right)\right)
\\\\ \ds~~~~~~
~=~(-1)\exp\left(\i{\pi\over 3}\right)~-~(-1)\exp\left(\i{5\pi\over 3}\right)~=~-\i\sqrt{3}.
\end{array}
\eeq
By Theorem One, we have
\bel{result2_x}
q_x(t,0)~=~-{\sqrt{3}\over 2\pi}\Gamma\left({1\over 3}\right)\int_0^t {\dot{g}(\tau)d\tau\over (t-\tau)^{1\over 3}}
\eeq
and
\bel{result2_xx}
 q_{xx}(t,0)~=~{\sqrt{3}\over 2\pi}\Gamma\left({2\over 3}\right)\int_0^t {\dot{g}(\tau)d\tau\over (t-\tau)^{2\over 3}}.
\eeq
Now, suppose $q_{xx}(t,0)=g(t)$ is given. We aim to find $q(t,0)$ and $q_x(t,0)$ in terms of $g(t)$. In this case, we have $u_1=2$, $v_1=0$ and $v_2=1$. 
\bel{V example2'}
\V(\rho)~=~\left[\begin{array}{lr}\ds
\rho^2\qquad\rho
\\\\ \ds
\rho\qquad \rho^2
\end{array}\right],\qquad 
\V^{11}(\rho)~=~\left[\begin{array}{lr}\ds
1\qquad\rho
\\\\ \ds
1\qquad \rho^2
\end{array}\right],\qquad 
\V^{21}(\rho)~=~\left[\begin{array}{lr}\ds
\rho^2\qquad1
\\\\ \ds
\rho\qquad 1
\end{array}\right].
\eeq
From computations, we have
\bel{quotients detV example2'}
{\det\V^{11}\over\det\V}(\rho)~=~ {\det\V^{21}\over\det\V}(\rho)~=~{1\over \rho(\rho+1)}~=~-1.
\eeq
By (\ref{Lambda_jl}), we have 
\bel{Lambda_11(2)'}
\begin{array}{lr}\ds
\Lambda_{11}~=~
\left({\det\V^{11}\over\det\V}\right)(\rho)\exp\left(\i(v_1-u_1)\left(\vartheta_1^1-{\pi\over 6}\right)\right)-\left({\det\V^{11}\over\det\V}\right)(\rho)\exp\left(\i(v_1-u_1)\left(\vartheta_1^2+{\pi\over 6}\right)\right)
\\\\ \ds~~~~~~
~=~(-1)\exp\left(\i{5\pi\over 3}\right)~-~(-1)\exp\left(\i{\pi\over 3}\right)~=~\sqrt{3}\i,
\end{array}
\eeq
and
\bel{Lambda_21(2)'}
\begin{array}{lr}\ds
\Lambda_{21}~=~
\left({\det\V^{21}\over\det\V}\right)(\rho)\exp\left(\i(v_2-u_1)\left(\vartheta_1^1-{\pi\over 6}\right)\right)-\left({\det\V^{21}\over\det\V}\right)(\rho)\exp\left(\i(v_2-u_1)\left(\vartheta_1^2+{\pi\over 6}\right)\right)
\\\\ \ds~~~~~~
~=~(-1)\exp\left(\i{11\pi\over 6}\right)~-~(-1)\exp\left(\i{7\pi\over 6}\right)~=~-\sqrt{3}.
\end{array}
\eeq
By Theorem One, we have
\bel{result2'}
q(t,0)~=~{\sqrt{3}\over 2\pi}\Gamma\left({1\over 3}\right)\int_0^t {g(\tau)d\tau\over (t-\tau)^{1\over 3}}
\eeq
and
\bel{result2'_x}
 q_{x}(t,0)~=~-{\sqrt{3}\over 2\pi}\Gamma\left({2\over 3}\right)\int_0^t {g(\tau)d\tau\over (t-\tau)^{2\over 3}}.
\eeq
Lastly, suppose $q_{x}(t,0)=g(t)$ is given. We aim to find $q(t,0)$ and $q_{xx}(t,0)$ in terms of $g(t)$. In this case, we have $u_1=1$, $v_1=0$ and $v_2=2$. 
\bel{V example2''}
\V(\rho)~=~\left[\begin{array}{lr}\ds
\rho^2\qquad1
\\\\ \ds
\rho\qquad 1
\end{array}\right],\qquad 
\V^{11}(\rho)~=~\left[\begin{array}{lr}\ds
\rho\qquad1
\\\\ \ds
\rho^2\qquad 1
\end{array}\right],\qquad 
\V^{21}(\rho)~=~\left[\begin{array}{lr}\ds
\rho^2\qquad\rho
\\\\ \ds
\rho\qquad \rho^2
\end{array}\right].
\eeq
From computations, we have
\bel{quotients detV example2''}
{\det\V^{11}\over\det\V}(\rho)~=~ {\det\V^{21}\over\det\V}(\rho)~=~\rho(\rho+1)~=~-1.
\eeq
By (\ref{Lambda_jl}), we have 
\bel{Lambda_11(2)''}
\begin{array}{lr}\ds
\Lambda_{11}~=~
\left({\det\V^{11}\over\det\V}\right)(\rho)\exp\left(\i(v_1-u_1)\left(\vartheta_1^1-{\pi\over 6}\right)\right)-\left({\det\V^{11}\over\det\V}\right)(\rho)\exp\left(\i(v_1-u_1)\left(\vartheta_1^2+{\pi\over 6}\right)\right)
\\\\ \ds~~~~~~
~=~(-1)\exp\left(\i{11\pi\over 6}\right)~-~(-1)\exp\left(\i{7\pi\over 6}\right)~=~-\sqrt{3},
\end{array}
\eeq
and
\bel{Lambda_21(2)''}
\begin{array}{lr}\ds
\Lambda_{21}~=~
\left({\det\V^{21}\over\det\V}\right)(\rho)\exp\left(\i(v_2-u_1)\left(\vartheta_1^1-{\pi\over 6}\right)\right)-\left({\det\V^{21}\over\det\V}\right)(\rho)\exp\left(\i(v_2-u_1)\left(\vartheta_1^2+{\pi\over 6}\right)\right)
\\\\ \ds~~~~~~
~=~(-1)\exp\left(\i{\pi\over 6}\right)~-~(-1)\exp\left(\i{5\pi\over 6}\right)~=~-\sqrt{3}.
\end{array}
\eeq
By Theorem One, we have
\bel{result2''}
q(t,0)~=~-{\sqrt{3}\over 2\pi}\Gamma\left({1\over 3}\right)\int_0^t {g(\tau)d\tau\over (t-\tau)^{1\over 3}}
\eeq
and
\bel{result2''_xx}
 q_{xx}(t,0)~=~-{\sqrt{3}\over 2\pi}\Gamma\left({2\over 3}\right)\int_0^t {\dot{g}(\tau)d\tau\over (t-\tau)^{2\over 3}}.
\eeq

{\bf Example 3:} Consider a second linear evolution equation with a third order derivative:
\bel{Third2}
q_t(t,x)-q_{xxx}(t,x)~=~0,\qquad0<x<\infty,\qquad0<t<T.
\eeq
Regarding to (\ref{Third2}), we have $a_n=\i$, $n=3$ and  $N=1$. The principal domain $\D$ is shown as below, where the rotation $\rho=e^{2\pi\i\over 3}$ maps $\D^+_2$ to $\D^-$ and $\rho^2$ maps $\D^+_1$ to $\D^-$.
 \begin{figure}[h]
\centering
\includegraphics[scale=0.25]{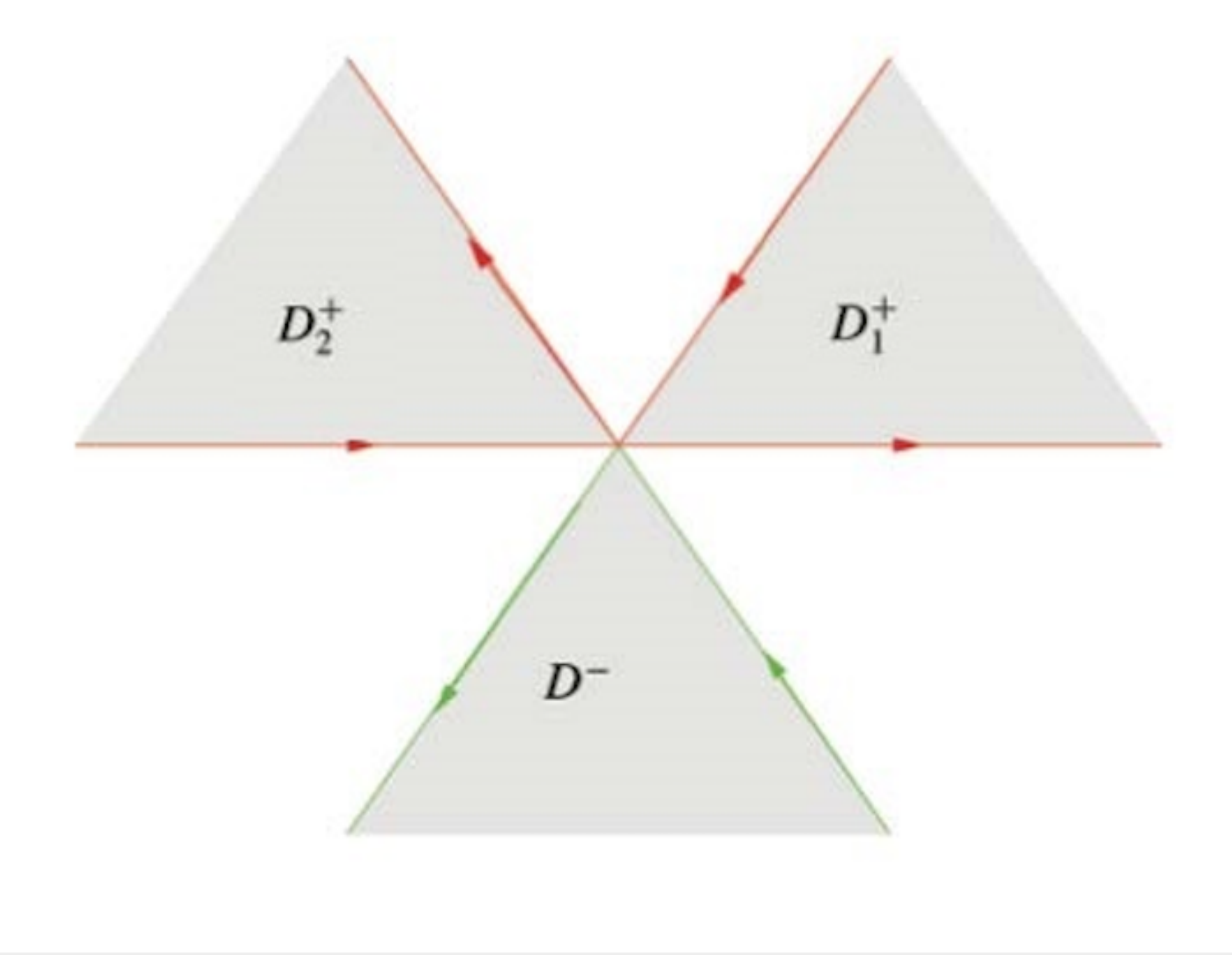}
\caption{$\p\D^+$ is represented as the red contour .}
\end{figure}

By taking $n=3$ and $\varphi={3\pi\over 2}$ in (\ref{vartheta Est}), we find 
\bel{vartheta(3)}
\vartheta_1^1~=~0\qquad\hbox{and}\qquad \vartheta_1^2~=~\pi.
\eeq
Let $q(t,0)=g_1(t)$ and $q_x(t,0)=g_2(t)$ be given. We aim to find $q_{xx}(t,0)$ in terms of $g_1(t)$ and $g_2(t)$. In this case, we have $u_1=0$, $u_2=1$ and $v_1=2$. Direct computations show that 
\bel{V example3}
\begin{array}{lr}\ds
\V(\rho)~=~1,\qquad \V(\rho^2)~=~1,
\\\\ \ds
\V^{11}(\rho)~=~\rho^2,\qquad \V^{11}(\rho^2)~=~\rho,
\\\\ \ds
\V^{12}(\rho)~=~\rho,\qquad \V^{12}(\rho^2)~=~\rho^2. 
\end{array}
\eeq
By (\ref{Lambda_jl}), we have 
\bel{Lambda_11(3)}
\begin{array}{lr}\ds
\Lambda_{11}~=~
\left({\det\V^{11}\over\det\V}\right)(\rho^2)\exp\left(\i(v_1-u_1)\left(\vartheta_1^1-{\pi\over 6}\right)\right)-\left({\det\V^{11}\over\det\V}\right)(\rho)\exp\left(\i(v_1-u_1)\left(\vartheta_1^2+{\pi\over 6}\right)\right)
\\\\ \ds~~~~~~
~=~\exp\left(\i{\pi\over 3}\right)~-~\exp\left(\i{5\pi\over 3}\right)~=~\i\sqrt{3},
\end{array}
\eeq
and
\bel{Lambda_12(3)}
\begin{array}{lr}\ds
\Lambda_{12}~=~
\left({\det\V^{12}\over\det\V}\right)(\rho^2)\exp\left(\i(v_2-u_1)\left(\vartheta_1^1-{\pi\over 6}\right)\right)-\left({\det\V^{12}\over\det\V}\right)(\rho)\exp\left(\i(v_2-u_1)\left(\vartheta_1^2+{\pi\over 6}\right)\right)
\\\\ \ds~~~~~~
~=~\exp\left(\i{\pi\over 6}\right)~-~\exp\left(\i{5\pi\over 6}\right)~=~-\sqrt{3}.
\end{array}
\eeq
By Theorem One, we have
\bel{result3_xx}
 q_{xx}(t,0)~=~-{\sqrt{3}\over 2\pi}\Gamma\left({2\over 3}\right)\int_0^t {\dot{g}_1(\tau)d\tau\over (t-\tau)^{2\over 3}}~-~{\sqrt{3}\over 2\pi}\Gamma\left({1\over 3}\right)\int_0^t {\dot{g}_2(\tau)d\tau\over (t-\tau)^{1\over 3}}.
\eeq

On the other hand, let $q_x(t,0)=g_1(t)$ and $q_{xx}(t,0)=g_2(t)$ be given. We aim to find $q(t,0)$ in terms of $g_1(t)$ and $g_2(t)$. In this case, we have $u_1=1$, $u_2=2$ and $v_1=0$. Direct computations show that 
\bel{V example3'}
\begin{array}{lr}\ds
\V(\rho)~=~\rho^2,\qquad \V(\rho^2)~=~\rho,
\\\\ \ds
\V^{11}(\rho)~=~\rho,\qquad \V^{11}(\rho^2)~=~\rho^2,
\\\\ \ds
\V^{12}(\rho)~=~1,\qquad \V^{12}(\rho^2)~=~1. 
\end{array}
\eeq
By (\ref{Lambda_jl}), we have 
\bel{Lambda_11(3)'}
\begin{array}{lr}\ds
\Lambda_{11}~=~
\left({\det\V^{11}\over\det\V}\right)(\rho^2)\exp\left(\i(v_1-u_1)\left(\vartheta_1^1-{\pi\over 6}\right)\right)-\left({\det\V^{11}\over\det\V}\right)(\rho)\exp\left(\i(v_1-u_1)\left(\vartheta_1^2+{\pi\over 6}\right)\right)
\\\\ \ds~~~~~~
~=~\exp\left(\i{5\pi\over 6}\right)~-~\exp\left(\i{\pi\over 6}\right)~=~-\sqrt{3},
\end{array}
\eeq
and
\bel{Lambda_12(3)'}
\begin{array}{lr}\ds
\Lambda_{12}~=~
\left({\det\V^{12}\over\det\V}\right)(\rho^2)\exp\left(\i(v_1-u_2)\left(\vartheta_1^1-{\pi\over 6}\right)\right)-\left({\det\V^{12}\over\det\V}\right)(\rho)\exp\left(\i(v_1-u_2)\left(\vartheta_1^2+{\pi\over 6}\right)\right)
\\\\ \ds~~~~~~
~=~\exp\left(\i{5\pi\over 3}\right)~-~\exp\left(\i{\pi\over 3}\right)~=~-\i\sqrt{3}.
\end{array}
\eeq
By Theorem One, we have
\bel{result3'}
 q(t,0)~=~-{\sqrt{3}\over 2\pi}\Gamma\left({2\over 3}\right)\int_0^t {g_1(\tau)d\tau\over (t-\tau)^{2\over 3}}~-~{\sqrt{3}\over 2\pi}\Gamma\left({1\over 3}\right)\int_0^t {g_2(\tau)d\tau\over (t-\tau)^{1\over 3}}.
\eeq
Lastly, let $q(t,0)=g_1(t)$ and $q_{xx}(t,0)=g_2(t)$ be given. We aim to find $q_x(t,0)$ in terms of $g_1(t)$ and $g_2(t)$. In this case, we have $u_1=0$, $u_2=2$ and $v_1=1$. Direct computations show that 
\bel{V example3''}
\begin{array}{lr}\ds
\V(\rho)~=~\rho,\qquad \V(\rho^2)~=~\rho^2,
\\\\ \ds
\V^{11}(\rho)~=~\rho^2,\qquad \V^{11}(\rho^2)~=~\rho,
\\\\ \ds
\V^{12}(\rho)~=~1,\qquad \V^{12}(\rho^2)~=~1. 
\end{array}
\eeq
By (\ref{Lambda_jl}), we have 
\bel{Lambda_11(3)''}
\begin{array}{lr}\ds
\Lambda_{11}~=~
\left({\det\V^{11}\over\det\V}\right)(\rho^2)\exp\left(\i(v_1-u_1)\left(\vartheta_1^1-{\pi\over 6}\right)\right)-\left({\det\V^{11}\over\det\V}\right)(\rho)\exp\left(\i(v_1-u_1)\left(\vartheta_1^2+{\pi\over 6}\right)\right)
\\\\ \ds~~~~~~
~=~\exp\left(\i{7\pi\over 6}\right)~-~\exp\left(\i{11\pi\over 6}\right)~=~-\sqrt{3},
\end{array}
\eeq
and
\bel{Lambda_12(3)''}
\begin{array}{lr}\ds
\Lambda_{12}~=~
\left({\det\V^{12}\over\det\V}\right)(\rho^2)\exp\left(\i(v_1-u_2)\left(\vartheta_1^1-{\pi\over 6}\right)\right)-\left({\det\V^{12}\over\det\V}\right)(\rho)\exp\left(\i(v_1-u_2)\left(\vartheta_1^2+{\pi\over 6}\right)\right)
\\\\ \ds~~~~~~
~=~\exp\left(\i{5\pi\over 6}\right)~-~\exp\left(\i{\pi\over 6}\right)~=~-\sqrt{3}.
\end{array}
\eeq
By Theorem One, we have
\bel{result3''_x}
 q_x(t,0)~=~-{\sqrt{3}\over 2\pi}\Gamma\left({1\over 3}\right)\int_0^t {\dot{g}_1(\tau)d\tau\over (t-\tau)^{1\over 3}}~-~{\sqrt{3}\over 2\pi}\Gamma\left({2\over 3}\right)\int_0^t {g_2(\tau)d\tau\over (t-\tau)^{2\over 3}}.
\eeq

{\bf Example 4:} We next consider an example which is an variance of first third order evolution equation discussed previously, known as the  first Stokes equation, with given canonical boundary condition:
\bel{Stoke1}
\left\{
\begin{array}{lr}\ds
q_t(t,x)+q_{xxx}(t,x)+q_x(t,x)~=~0,\qquad0<x<\infty,\qquad0<t<T;
\\\\ \ds ~~~~~~~~~~~~~~~~~~~~~~~~~~~~~~~~~~~
q(t,0)~=~g(t),\qquad 0<t<T,
\end{array}\right.
\eeq
Regarding to (\ref{Stoke1}), we have $n=3$, $N=2$ and
\bel{omega Stoke1}
\omega(\xi)~=~-\i\xi^3+\i\xi.
\eeq
The domain $D$ regarding $\omega(\xi)$ in (\ref{omega Stoke1}) is shown as below:
  \begin{figure}[h]
\centering
\includegraphics[scale=0.40]{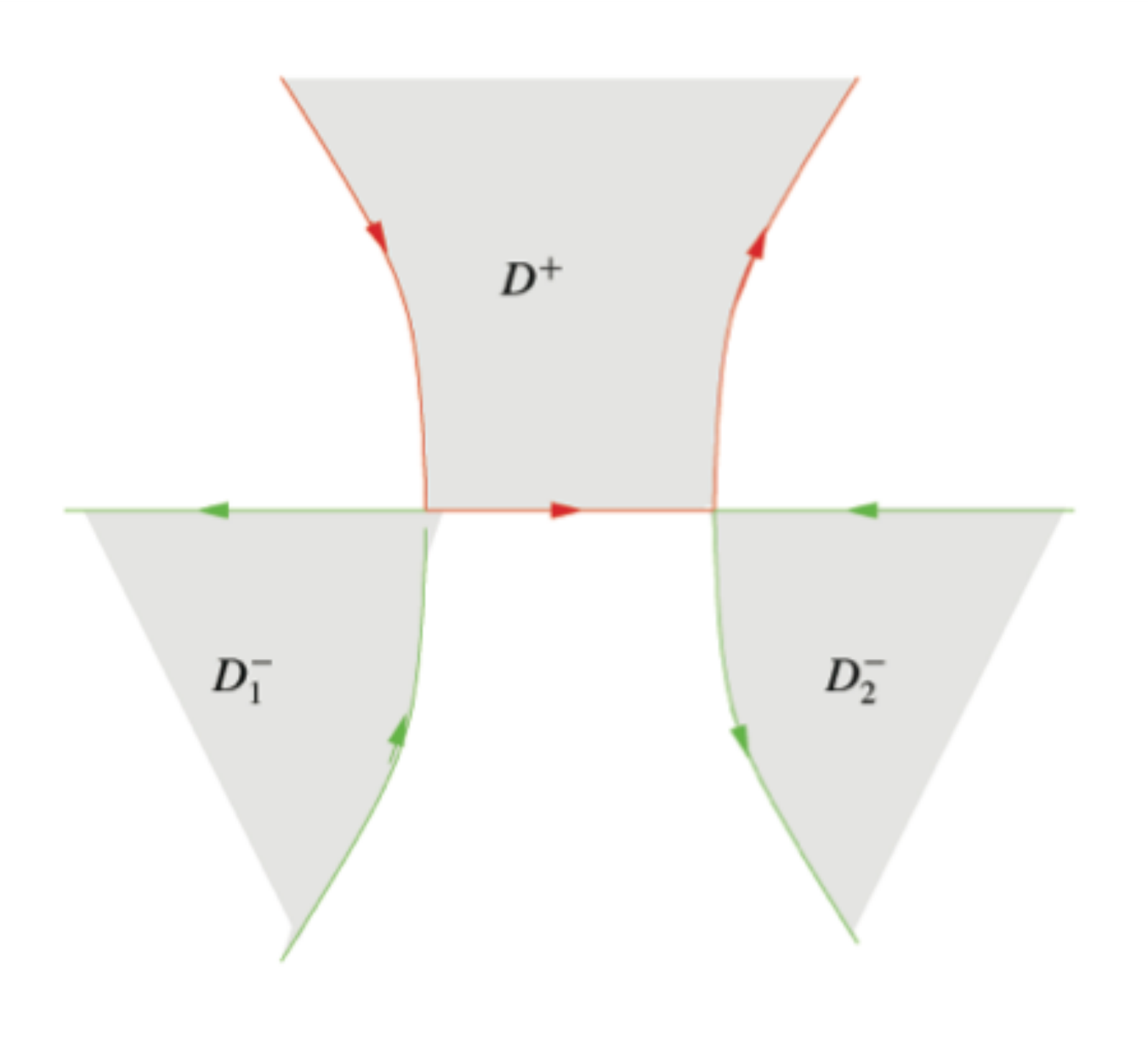}
\caption{$\p D^+$ is represented as the red contour .}
\end{figure}

We aim to find $q_x(t,0)$ and $q_{xx}(t,0)$ in terms of $g(t)$. Let $u_1=0$, $v_1=1$ and $v_2=2$. For $\xi\in\Cx$ fixed, the roots of (\ref{Poly}) can be determined by
\bel{Poly Stoke1}
{\omega(z)-\omega(\xi)\over z-\xi}~=~-\i\left(z^2+\xi z+\xi^2-1\right)~=~0.
\eeq
In this case, $D^+_L=D^+_{L,1}$. We abbreviate $z_i=z^1_i,~i=1,2$ for which
\bel{z_1,z_2 Stoke1}
z_1(\xi)~=~{-\xi+\left(4-3\xi^2\right)^{1\over 2}\over 2},
\qquad
z_2(\xi)~=~{-\xi-\left(4-3\xi^2\right)^{1\over 2}\over 2}.
\eeq
Notice that as $|\xi|\mt\infty$ for $\xi\in D^+$, we have $\Im z_1(\xi)\leq0$ and $\Im z_2(\xi)\leq0$. Recall from (\ref{omega_m}) and (\ref{A}), we have
\bel{A Stoke1}
\A(z_1,z_2)~=~\left[\begin{array}{lr}\ds
-\i z_1\quad-\i
\\\\ \ds
-\i z_2\quad -\i
\end{array}\right],\quad 
\A^{11}(z_1,z_2)~=~\left[\begin{array}{lr}\ds
-\i z_1^2+\i \quad-\i
\\\\ \ds
-\i z_2^2+\i \quad -\i
\end{array}\right],\quad 
\A^{21}(z_1,z_2)~=~\left[\begin{array}{lr}\ds
-\i z_1\quad -\i z_1^2+\i
\\\\ \ds
-\i z_2 \quad -\i z_2^2+\i
\end{array}\right].
\eeq
Direct computations show that 
\bel{quotients Stoke1}
{\det\A^{11}\over \det\A}(z_1,z_2)~=~z_1+z_2,\qquad
{\det\A^{21}\over \det\A}(z_1,z_2)~=~-\left(1+z_1z_2\right).
\eeq
From (\ref{z_1,z_2 Stoke1}) we find $z_1+z_2=-\xi$ and $1-z_1z_2=-\xi^2$ which are both analytic. 
$\omega(\xi)$ has zeros at $-1$, $0$ and $1$. By Theorem Two, we have
\bel{result Stoke1}
2\pi q_x(t,0)~=~-~\hbox{\bf p.v} \int_{-\infty}^\infty {3\xi^3-\xi\over\xi^3-\xi}\left(\int_0^t e^{\i(\xi^3-\xi)(t-\tau)} \dot{g}(\tau) d\tau\right) d\xi,
\eeq
\bel{result Stoke1'}
2\pi q_{xx}(t,0)~=~-2\pi g(t)~-~\hbox{\bf p.v} \int_{-\infty}^\infty \i{3\xi^4-\xi^2\over\xi^3-\xi}\left(\int_0^te^{\i(\xi^3-\xi)(t-\tau)} \dot{g}(\tau)d \tau\right) d\xi.
\eeq

{\bf Example 5:} Our last example is the second Stokes equation with given canonical boundary condition:  
\bel{Stoke2}
\left\{
\begin{array}{lr}\ds
q_t(t,x)-q_{xxx}(t,x)+q_x(t,x)~=~0,\qquad0<x<\infty,\qquad0<t<T;
\\\\ \ds ~~~~~~~~~~~~~~~~~~~~~~~~~~~~~~~~~~~
q(t,0)~=~g_1(t),\qquad 0<t<T,
\\\\ \ds ~~~~~~~~~~~~~~~~~~~~~~~~~~~~~~~~~~~
q_x(t,0)~=~g_2(t),\qquad 0<t<T.
\end{array}\right.
\eeq
Regarding to (\ref{Stoke2}), we have $n=3$, $N=1$ and
\bel{omega Stoke2}
\omega(\xi)~=~\i\xi^3+\i\xi.
\eeq
The domain $D$ regarding $\omega(\xi)$ in (\ref{omega Stoke2}) is shown as below:
  \begin{figure}[h]
\centering
\includegraphics[scale=0.35]{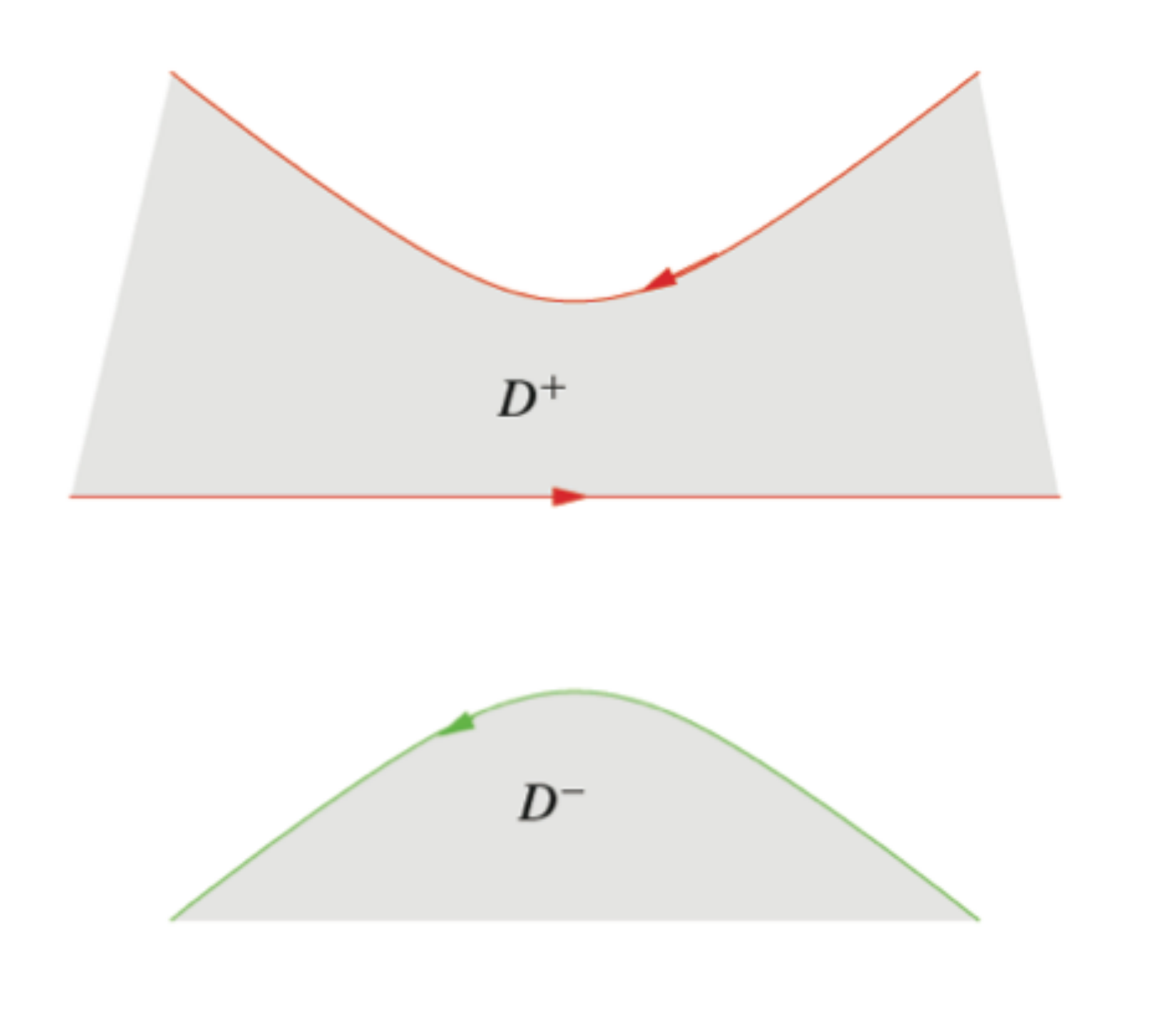}
\caption{$\p D^+$ is represented as the red contour .}
\end{figure}

Let $u_1=0$, $u_2=1$ and $v_1=2$. For $\xi\in\Cx$ fixed, the roots of (\ref{Poly}) can be determined by
\bel{Poly Stoke2}
{\omega(z)-\omega(\xi)\over z-\xi}~=~\i\left(z^2+\xi z+\xi^2+1\right)~=~0.
\eeq
We find that
\bel{z_1,z_2 Stoke2}
z_1(\xi)~=~{-\xi+\left(-4-3\xi^2\right)^{1\over 2}\over 2},
\qquad
z_2(\xi)~=~{-\xi-\left(-4-3\xi^2\right)^{1\over 2}\over 2}
\eeq
where  $\Im z_1(\xi)\leq0$ as $|\xi|\mt\infty$ for $\xi\in D^+_{L,1}$ and $\Im z_1^2(\xi)\leq0$ as $|\xi|\mt\infty$ for $\xi\in\D^+_{L,2}$. 

In this case, $D^-_L=D^-_{L,1}$. We abbreviate $z^k=z_1^k,~k=1,2$.
Recall from (\ref{omega_m}) and (\ref{A}). We have 
\bel{A Stoke2}
\begin{array}{lr}
\A(z^k)~=~\i,
\qquad
\A^{11}(z^k)~=~\i (z^k)^2+\i,
\qquad 
\A^{12}(z^k)~=~\i z^k
\\\\ \ds
{\det\A^{11}\over \det\A}(z^k)~=~1+(z^k)^2,\qquad
{\det\A^{12}\over \det\A}(z^k)~=~z^k,\qquad k=1,2.
\end{array}
\eeq
The quotients of determinants in (\ref{A Stoke2}) are analytic for $|\xi|>2/\sqrt{3}>1$, for which the discriminant in (\ref{z_1,z_2 Stoke2}) is nonvanishing. 
On the other hand, $\omega(\xi)=0$ in (\ref{omega Stoke2}) at $-\i$, $0$ and $\i$. 
By Theorem Two, we have
\bel{result Stoke2}
\begin{array}{lr}\ds
4\pi q_{xx}(t,0)~=~\pi g_2(t)~+~\hbox{\bf p.v}\int_{0}^{\infty} \Big(z^1(\xi)-z^1(\i\xi)+z^2(\i\xi)-z^2(\xi)\Big) {3\xi^2+1\over\xi^3+\xi}\left(\int_0^te^{\i(\xi^3+\xi)(t-\tau)} \dot{g}_2(\tau) d\tau\right)d\xi
\\\\ \ds ~~~~~~~~~~~~~~~~
~+~\hbox{\bf p.v} \int_{0}^{\infty} \i\Big((z^1)^2(\xi)-(z^1)^2(\i\xi)+(z^2)^2(\i\xi)-(z^2)^2(\xi)\Big) {3\xi^2+1\over\xi^3+\xi}\left(\int_0^te^{\i(\xi^3+\xi)(t-\tau)} \dot{g}_1(\tau) d\tau\right)d\xi.
\end{array}
\eeq

\end{document}